\documentclass[11pt,reqno]{amsart}
\usepackage{amsmath,amssymb,amsthm}
\usepackage[margin=1in]{geometry}
\usepackage{booktabs,array}
\usepackage{tikz}
\usetikzlibrary{positioning}
\usepackage[colorlinks=true,linkcolor=blue,citecolor=blue,urlcolor=blue]{hyperref}

\newtheorem{theorem}{Theorem}[section]
\newtheorem{proposition}[theorem]{Proposition}
\newtheorem{lemma}[theorem]{Lemma}
\newtheorem{corollary}[theorem]{Corollary}
\newtheorem{fact}[theorem]{Fact}
\newtheorem{mainthm}{Theorem}

\theoremstyle{definition}
\newtheorem{setting}[theorem]{Setting}
\theoremstyle{remark}
\newtheorem{remark}[theorem]{Remark}
\numberwithin{equation}{section}
\numberwithin{table}{section}
\numberwithin{figure}{section}
\newcommand{\Z}{\mathbb{Z}}
\newcommand{\R}{\mathbb{R}}
\newcommand{\C}{\mathbb{C}}
\newcommand{\Q}{\mathbb{Q}}
\newcommand{\e}{\mathrm{e}}
\newcommand{\eps}{\varepsilon}
\newcommand{\Ga}{\Gamma}
\newcommand{\hb}{\overline{h}}
\newcommand{\Db}{\overline{\Delta}}
\newcommand{\GQ}{G_Q}
\newcommand{\zQ}{\zeta_{|Q|}}
\newcommand{\Hs}{\mathbf H}
\newcommand{\vs}{\varsigma}
\newcommand{\Pc}{\mathcal P}
\newcommand{\Dc}{\mathcal D}
\DeclareMathOperator{\sgn}{sgn}
\DeclareMathOperator{\arsinh}{arsinh}
\DeclareMathOperator{\Var}{Var}
\DeclareMathOperator{\supp}{supp}

\begin{document}

\title[Gaps between zeros of Epstein zeta-functions]{Gaps between zeros of Epstein zeta-functions\\ of binary quadratic forms on the critical line}

\author{Stephan Baier}
\address{Stephan Baier, Ramakrishna Mission Vivekananda Educational Research Institute, Department of Mathematics,
G.~T.~Road, PO-Belur Math, Howrah, West Bengal 711202, India}
\email{stephanbaier2017@gmail.com}

\author{Srinivas Kotyada}
\address{Srinivas Kotyada, Indian Institute of Science Education and Research (IISER) Tirupati, Tirupati 517619,
Andhra Pradesh, India}
\email{ksrinivas@labs.iisertirupati.ac.in}

\author{Usha Keshav Sangale}
\address{Usha Keshav Sangale, SRTM University, Nanded, Maharashtra 431606, India}
\email{ushas073@gmail.com}

\subjclass[2020]{Primary 11E45, 11M41; Secondary 11L07, 11R42}
\keywords{Epstein zeta-function, zeros on the critical line, gaps between zeros, Voronoi summation, indefinite
quadratic forms, real quadratic fields, two-dimensional exponential sums}

\begin{abstract}
We prove that for every integral binary quadratic form $Q$ whose discriminant is, up to sign, not a perfect square,
the Epstein zeta-function of $Q$ has a zero $\frac12+i\gamma$ with $|\gamma-T|\le T^{2/5+\eps}$ for all large $T$. For
indefinite $Q$, that is, for zeta-functions of ideal classes of real quadratic orders, this is the first exponent
below the classical barrier $1/2$. For positive definite $Q$ it improves the exponents $5/11$ of Jutila and Srinivas
and $3/7$ of the authors, and for discriminant $-4m^2$ with $m$ odd we obtain $1/3$. A main ingredient, of independent
interest, is an explicit Voronoi summation formula with additive twists to every modulus for the representation
numbers of an indefinite binary form. The resulting dual sums are treated by a two-dimensional Weyl differencing,
controlled by an exact formula for the Hessian determinant of the differenced phase.
\end{abstract}

\maketitle

\section{Introduction}\label{sec:intro}

\subsection{Epstein zeta-functions and their zeros}
Let $Q(x,y)=ax^2+bxy+cy^2$ be a positive definite integral binary quadratic form, $\Delta=4ac-b^2$, and let
\[
\zeta_Q(s)=\sum_{(m,n)\in\Z^2\setminus\{0\}}Q(m,n)^{-s}=\sum_{n\ge1}r_Q(n)n^{-s}\qquad(\sigma>1)
\]
be its Epstein zeta-function \cite{Ep}. It continues meromorphically to $\C$, with a simple pole at $s=1$, and the
completed function $(\sqrt\Delta/2\pi)^s\Ga(s)\zeta_Q(s)$ is invariant under $s\mapsto1-s$. Through the correspondence
between forms and ideal classes, $\zeta_Q$ is (up to a constant) the zeta-function of an ideal class of an imaginary
quadratic order, and hence a linear combination of Hecke $L$-functions of class group characters. It has an Euler
product only in exceptional cases. When the class number exceeds one, Davenport and Heilbronn \cite{DH} showed that
$\zeta_Q$ has infinitely many zeros in $\sigma>1$, so there is no Riemann hypothesis for $\zeta_Q$. Nevertheless, Potter
and Titchmarsh \cite{PT} proved that $\zeta_Q$ has infinitely many zeros on the critical line, and Bombieri and Hejhal
\cite{BH} showed, under the generalized Riemann hypothesis and a spacing hypothesis, that almost all zeros of such
linear combinations lie on it.

This paper is concerned with the \emph{gaps} between consecutive zeros on the critical line. For the Riemann
zeta-function, Hardy and Littlewood \cite{HL} proved that every interval $[T,T+T^{1/4+\eps}]$ contains the ordinate
of a zero on the critical line; the exponent has since been reduced below $1/6$ (Karatsuba \cite{Ka}; see also
\cite{Iv}). For Epstein zeta-functions the analogue of the Hardy--Littlewood exponent is $1/2$: Potter and Titchmarsh
showed in 1935 that $\zeta_Q$ has a zero $\frac12+i\gamma$ with $\gamma\in[T,T+T^{1/2+\eps}]$, and Sankaranarayanan \cite{Sa}
replaced $T^\eps$ by $c\log T$. For seventy years no smaller exponent was known. Jutila and Srinivas \cite{JS} were the
first to go below $1/2$, with the exponent $5/11$, and the authors \cite{BSS} improved it to $3/7$. For holomorphic cusp
forms of level one, Jutila \cite{J1} had obtained the exponent $1/3$; we explain below why that argument does not
extend to Epstein zeta-functions.

For an \emph{indefinite} form the Epstein zeta-function is defined modulo the infinite group of automorphs, and it is
the zeta-function of an ideal class of a real quadratic order. A Potter--Titchmarsh type result for Dedekind
zeta-functions of real quadratic fields goes back to Berndt \cite{B}, and Sankaranarayanan's theorem \cite{Sa} gives the
exponent $1/2$ for ideal class zeta-functions of real quadratic fields. In the concluding remark of \cite{JS} the question
was raised whether the exponent $1/2$ can be broken in the real quadratic case as well; the method of \cite{JS} does
not apply there, because it rests on a Voronoi formula for the representation numbers of a positive definite form. We
answer this question (Theorem~\ref{thm:B}), and at the same time improve the exponent in the definite case
(Theorem~\ref{thm:A}). The new Voronoi formula (Theorem~\ref{thm:voronoi}) on which the indefinite case rests seems to
be of independent interest.

\subsection{Why the exponent $1/2$ is a barrier}\label{subsec:difficulty}
All results below $1/2$ rest on the method of Hardy and Littlewood. If the Hardy function $W(t)$ of $\zeta_Q$ (real on
the real line) had no zero in $[T-H,T+H]$, then $|\int W(T+u)\e^{-(u/H_0)^2}du|$ would equal $\int|W(T+u)|\e^{-(u/H_0)^2}du$,
which is large by a lower bound of Ramachandra. An approximate formula for $\zeta_Q$ and Stirling's formula show that
only the terms near the \emph{resonance point} $x_*=T\sqrt\Delta/2\pi$ contribute to the integral. So it suffices to
show that the \emph{critical Dirichlet polynomial}
\begin{equation}\label{eq:crit}
\Pc(t)=\sum_n\eta(n)r(n)n^{-1/2-it},\qquad\supp\eta\subset[x_*-K,x_*+K],\qquad HK=T^{1+2\eps},
\end{equation}
is $o(1)$ for $|t-T|\le H$. The trivial bound $\Pc(t)\ll KT^{-1/2}$ gives $H=T^{1/2+\eps}$, i.e.\ the exponent $1/2$. Any
smaller exponent requires a power saving in a sum of length $K>T^{1/2}$ whose terms oscillate like $n^{-it}$, $t\asymp T$.

\emph{The arithmetic obstruction.} To transform $\Pc(t)$ one twists $n^{-it}=\e(-hn/k)\cdot\e(hn/k)n^{-it}$ with a
fraction $h/k$ for which $\e(hx/k)x^{-it}$ is stationary inside the support of $\eta$; this forces $h/k$ to be within
$O(K/T)$ of the \emph{critical ratio} $1/\sqrt\Delta$. A Voronoi formula with additive twist, followed by stationary
phase, then turns $\Pc(t)$ into a dual sum of length
\[
N_0\asymp\frac{h^2K^2}T .
\]
For level one cusp forms the critical ratio is $1$, one takes $h/k=1/1$, and $N_0\asymp K^2/T$ is short; the trivial bound
for the dual sum already gives the exponent $1/3$ \cite{J1}. For Epstein zeta-functions with $\Delta$ not a square, the
critical ratio $1/\sqrt\Delta$ is a quadratic irrational. By Lagrange's theorem $|1/\sqrt\Delta-h/k|\gg k^{-2}$, hence
$k\gg\sqrt{T/K}$ and
\[
N_0\asymp K .
\]
\emph{The dual sum is exactly as long as the original one.} The transformation maps $\Pc(t)$ to a two-dimensional
exponential sum $S$ \emph{at its self-dual point}, with $\asymp K$ terms and a phase of size $\asymp K$; at this point the
Voronoi step, like the B-process of van der Corput, is an involution and yields no saving by itself. Every saving must
therefore come from genuine cancellation in $S$, and the exponent is governed by the saving $\sigma$ in a bound
$|S|\ll K^{1-\sigma}$:
\begin{equation}\label{eq:theta}
\theta=\frac{1-2\sigma}{2(1-\sigma)} .
\end{equation}
\begin{center}\small
\begin{tabular}{@{}lccl@{}}
\toprule
estimate of the dual sum & $\sigma$ & $\theta$ & source\\
\midrule
trivial & $0$ & $1/2$ & \cite{PT}, \cite{Sa}\\
third derivative test in one variable & $1/12$ & $5/11$ & \cite{JS}\\
Weyl differencing in one variable, B-process & $1/8$ & $3/7$ & \cite{BSS}\\
Weyl differencing in two variables & $1/6$ & $2/5$ & this paper\\
\bottomrule
\end{tabular}
\end{center}

\emph{The analytic obstruction.} The phase of $S$ is $\varpi P(w)+\vs\Phi(P(w))$, where $P$ is the adjoint form. The first
term is a quadratic polynomial whose leading coefficient $\varpi$ is essentially a rational number with denominator
$k\asymp\sqrt{T/K}$, so it can only be removed by differencing. The second term is, up to a negligible correction, $2\kappa\sqrt{P(w)}$
with $\kappa\asymp K^{1/2}$; it is homogeneous of degree one, so its Hessian has rank one and degenerates in the radial
direction, and the standard two-dimensional van der Corput estimates do not apply.

\emph{The limit of one-variable methods.} The methods of \cite{JS} and \cite{BSS} estimate the sum over one variable
non-trivially. In the scheme of \cite{BSS} the Weyl diagonal costs a factor which no estimate in the other variable can
recover: even the conjectural exponent pair $(0,\frac12)$ would only give $\theta=2/5$ (Remark~\ref{rem:limits}). Thus
the exponent $2/5$ of Theorems~\ref{thm:A} and~\ref{thm:B} is what the one-variable approach would give under the exponent
pair conjecture; we obtain it unconditionally, by treating both variables at once. Going beyond $2/5$ requires a
saving beyond Weyl strength for a two-dimensional sum at its self-dual point.

\emph{The indefinite case.} Here the obstacles appear already before the exponential sum. The representation numbers
count orbits under an infinite group, the theta series of $Q$ is not holomorphic (by Maass \cite{Ma}, the relevant
$L$-functions belong to Maass forms of weight $0$), the zeta-function of a narrow ideal class has no self-dual
functional equation (Remark~\ref{rem:narrow}), and no explicit Voronoi formula with additive twists was available.

\subsection{Results}
Let now $Q(x,y)=ax^2+bxy+cy^2\in\Z[x,y]$ be indefinite with $\Delta=b^2-4ac>0$ not a square. Then $Q(v)\ne0$ for
$v\in\Z^2\setminus\{0\}$, and every represented integer is represented infinitely often, because of the group
$\Ga=SO_Q(\Z)=\{\gamma\in SL_2(\Z):Q\circ\gamma=Q\}$ of proper automorphs. Put
\[
r_Q(n)=\#\bigl(\Ga\backslash\{v\in\Z^2:Q(v)=n\}\bigr)\qquad(n\in\Z\setminus\{0\})
\]
and, for $\sigma>1$,
\[
\zQ(s)=\sum_{n\ge1}\frac{r_Q(n)+r_Q(-n)}{n^s},\qquad\Xi^-_Q(s)=\sum_{n\ge1}\frac{r_Q(n)-r_Q(-n)}{n^{s}} .
\]
If $Q$ is primitive and $\Delta$ is fundamental, $\zQ$ is the zeta-function of a wide ideal class of $\Q(\sqrt\Delta)$ if every
unit has norm $+1$, and twice it otherwise (Remark~\ref{rem:ideals}). It continues to $\C$ with a single simple pole at
$s=1$, and
\begin{equation}\label{eq:FEintro}
\Lambda(s):=\Bigl(\frac{\sqrt\Delta}{\pi}\Bigr)^s\Ga\Bigl(\frac s2\Bigr)^2\zQ(s)=\Lambda(1-s)
\end{equation}
(Corollary~\ref{cor:selfdual}); $\Xi^-_Q$ is entire and satisfies the analogous functional equation with
$\Ga(\frac{s+1}2)^2$. Our main results on zeros are the following.

\begin{mainthm}[Positive definite forms]\label{thm:A}
Let $Q$ be a positive definite integral binary quadratic form with $\Delta=4ac-b^2$, and let $\eps>0$.
\begin{enumerate}
\item[(a)] If $\Delta$ is not a perfect square, then for $T\ge T_0(\eps,Q)$ the function $\zeta_Q$ has a zero
$\frac12+i\gamma$ with $|\gamma-T|\le T^{2/5+\eps}$.
\item[(b)] If $Q$ is primitive and $\Delta=4m^2$ with $m$ odd, then for $T\ge T_0(\eps,Q)$ the function $\zeta_Q$ has a zero
$\frac12+i\gamma$ with $|\gamma-T|\le T^{1/3+\eps}$.
\end{enumerate}
\end{mainthm}

\begin{mainthm}[Indefinite forms]\label{thm:B}
Let $Q$ be an indefinite integral binary quadratic form whose discriminant $\Delta=b^2-4ac$ is not a perfect square,
and let $\eps>0$. Then for $T\ge T_0(\eps,Q)$ the function $\zQ$ has a zero $\frac12+i\gamma$ with $|\gamma-T|\le T^{2/5+\eps}$.
The same holds for $\Xi^-_Q$, provided $\Xi^-_Q\not\equiv0$.
\end{mainthm}

The key to Theorem~\ref{thm:B} is a Voronoi formula for $r_Q(n)$ twisted by $\e(hn/k)$. Let
\[
Q^*(w)=cw_1^2-bw_1w_2+aw_2^2
\]
be the adjoint form, $\Ga^*$ its group of proper automorphs, $\GQ(k,h,w)=\sum_{x\in(\Z/k\Z)^2}\e_k(hQ(x)+x\cdot w)$ the
Gauss sum, $\GQ(k,h)=\GQ(k,h,0)$, and
\[
\rho(n)=\rho_Q(n;h/k)=\sum_{w\in\Ga^*\backslash\Z^2,\ Q^*(w)=n}\GQ(k,h,w)\qquad(n\in\Z\setminus\{0\});
\]
the sum is well defined by Lemma~\ref{lem:gauss}(a). Let $\eps_0>1$ be the fundamental unit of norm $1$ attached to $Q$
(\S\ref{subsec:automorphs}).

\begin{mainthm}[Voronoi formula for indefinite forms]\label{thm:voronoi}
Let $Q$ be indefinite with non-square discriminant $\Delta$, let $f\in C^3$ have compact support in $(0,\infty)$, and
let $k\ge1$, $(h,k)=1$. Then
\begin{multline*}
\sum_{n\ge1}r_Q(n)\e\Bigl(\frac{hn}{k}\Bigr)f(n)=\frac{\log\eps_0\,\GQ(k,h)}{k^2\sqrt\Delta}\int_0^\infty f(x)\,dx\\
+\frac1{k^2\sqrt\Delta}\sum_{n\ge1}\Bigl\{4\rho(n)\int_0^\infty f(x)K_0\Bigl(\frac{4\pi\sqrt{nx}}{k\sqrt\Delta}\Bigr)dx
-2\pi\rho(-n)\int_0^\infty f(x)Y_0\Bigl(\frac{4\pi\sqrt{nx}}{k\sqrt\Delta}\Bigr)dx\Bigr\}.
\end{multline*}
\end{mainthm}

Some comments are in order.

(i) Theorem~\ref{thm:B} applies to the zeta-function of every wide ideal class of every real quadratic order. To our
knowledge it is the first result with an exponent below $1/2$ in the real quadratic case.

(ii) Theorem~\ref{thm:voronoi} is the indefinite counterpart of the classical Voronoi formula for positive definite forms
(Proposition~\ref{prop:Vdef}), whose kernel is $J_0$. For indefinite forms the kernels are exactly those of Voronoi's
formula for the divisor problem, $-2\pi Y_0$ and $4K_0$: the representations of negative integers by $Q^*$ feed the
oscillatory kernel, and those of positive integers the exponentially decaying one. We emphasise three features. The
formula holds for \emph{every} modulus $k$, including $(k,\Delta)>1$; the dual coefficients are explicit arithmetic
quantities, Gauss sums of $Q$ summed over the $\Ga^*$-orbits of the adjoint form; and the proof (Theorem~\ref{thm:FE})
is elementary, by Hecke's method applied to shifted lattices with sign characters. Abstractly, the theta series of $Q$
is a combination of Maass forms attached to Hecke characters of a real quadratic field \cite{Ma}, and Voronoi formulas
for $\mathrm{GL}(2)$ newforms with arbitrary additive twists are known \cite{KZ}; but deducing Theorem~\ref{thm:voronoi}
in this way would require the explicit decomposition into newforms, including old forms and imprimitive characters,
which Theorem~\ref{thm:voronoi} avoids. We expect the formula to be useful for other problems on indefinite forms
and real quadratic fields, such as the lattice point problem for $\sum_{n\le x}r_Q(n)$ in hyperbolic regions, the
distribution of $r_Q(n)$ in arithmetic progressions and short intervals, and additive problems involving $r_Q(n)$.

(iii) If $\Delta$ is a square, then $\Delta=4m^2$ (as $\Delta\equiv0,3\bmod4$), and these are exactly the forms attached to orders in
$\Z[i]$. The critical ratio $1/\sqrt\Delta=1/(2m)$ is then rational, and Theorem~\ref{thm:A}(b) is the analogue of Jutila's
exponent for cusp forms. For \emph{even} $m$, the main term of the transformed Dirichlet polynomial produces a non-zero
constant in the Hardy function, of the same size as Ramachandra's lower bound, so the method cannot work
(Proposition~\ref{prop:square}). For $Q=x^2+4y^2$ the zeros are nevertheless spaced $O(1)$ apart, for a different reason
(Proposition~\ref{prop:x24y2}).

(iv) For a narrow ideal class $C^+$ of a real quadratic field with $C^+\ne C^+[(\sqrt\Delta)]$, the zeta-function
$\zeta(s,C^+)=\frac12(\zQ+\Xi^-_Q)$ is not a multiple of a self-dual function, and Hardy's method does not apply to it
(Remark~\ref{rem:narrow}). Theorem~\ref{thm:B} covers its even and odd parts.

(v) In Remark~\ref{rem:BSS} we correct an error term in \cite[(4.9)]{BSS}. The exponent $3/7$ of \cite{BSS} is not affected.

\subsection{The method}\label{subsec:method}
The proof follows the scheme just described: zero detection, transformation of $\Pc(t)$, estimation of a double
exponential sum. Section~\ref{sec:strategy} gives a complete outline with all key formulas; here we point out what is new.

\emph{Twisted functional equations and Voronoi formulas for every modulus.} In \cite{JS} the transformation of $\Pc(t)$
is stated as a direct analogue of Jutila's transformation formula for the divisor function \cite[Thm.~4.3]{J2} (see also \cite{J3}), which
rests on Voronoi's formula and a specific smoothing $\eta_J$. We give a self-contained proof, uniform in the definite and
indefinite cases. For indefinite $Q$ we prove, by Hecke's method applied to shifted lattices with sign characters, a
functional equation for $\sum_nr_Q(n)\e(hn/k)n^{-s}$ for \emph{every} modulus $k$ (Theorem~\ref{thm:FE}), which leads
to Theorem~\ref{thm:voronoi}. In the transformation of $\Pc(t)$ only the $Y_0$-part survives the saddle point analysis. The integral transforms are then evaluated by a stationary phase lemma
of Blomer, Khan and Young \cite{BKY} (Lemma~\ref{lem:statphase} and Proposition~\ref{prop:transf}).

\emph{Two-dimensional Weyl differencing.} The step that lowers the exponent is to difference the dual sum in both
variables, with shifts $\ell\in\Z^2$, $|\ell|_\infty<M$ (Figure~\ref{fig:shifts}). The diagonal loss becomes $N^2/M^2$ instead
of $N^2/M$, while each differenced sum is still bounded essentially by $\kappa|\ell|$, where $\kappa\asymp K^{1/2}$. The
difficulty is that the differenced phase $g_\ell(w)\approx2\kappa(\sqrt{P(w+\ell)}-\sqrt{P(w)})$ has a Hessian which degenerates
when $w$ is parallel to $\ell$. We control it through the exact identity (Lemma~\ref{lem:hess})
\[
\det\nabla^2\Bigl(\sqrt{P(w+\ell)}-\sqrt{P(w)}\Bigr)=-\frac{d^2}{16}\,\frac{\det(\ell,w)^2}{P(w+\ell)^{3/2}P(w)^{3/2}},\qquad d=\operatorname{disc}P,
\]
together with a decomposition of the summation region according to the angle between $w$ and $\ell$. On each piece we
apply the B-process along a direction transverse to $w$ and van der Corput's estimate in the other variable. The number
of intervals that arise is bounded uniformly by o-minimality \cite{vdD,vdDM}. The resulting bound
$S\ll N^{2/3}K^{1/6}\log T$ (Theorem~\ref{thm:SN}) is a saving of Weyl strength, $\sigma=\frac16$.

\emph{The square case.} For $\Delta=4m^2$ we may take $h/k=1/(2m)$, so that the dual sum is short, $N_0\asymp K^2/T$. The
main term of the transformation is a multiple of the Gauss sum $\GQ(2m,-1)$, which vanishes exactly when $m$ is odd
(Lemma~\ref{lem:Gsq}).

\emph{Limits.} The exponent $1/3$ for non-square $\Delta$ would require the saving $\sigma=\frac14$ for a two-dimensional sum at
its self-dual point. This is beyond Weyl strength; for $\mathrm{GL}(2)$ in the $t$-aspect the Weyl barrier has only recently been
crossed \cite{HMSS,HMQ}. See Remark~\ref{rem:limits}.

\subsection{Organisation and conventions}
Section~\ref{sec:strategy} outlines the proof. Section~\ref{sec:forms} collects facts on binary forms, automorphs and Gauss
sums, and Section~\ref{sec:FEV} proves the twisted functional equations and Voronoi formulas. Section~\ref{sec:osc} contains
the stationary phase analysis, and Section~\ref{sec:transf} the zero detection and the transformation of $\Pc(t)$.
Sections~\ref{sec:lemmas} and \ref{sec:expsum} estimate the double exponential sum, and Section~\ref{sec:proof} completes the
proofs of Theorem~\ref{thm:A}(a) and Theorem~\ref{thm:B}. Section~\ref{sec:square} treats square discriminants, and
Section~\ref{sec:remarks} contains further remarks.

$\e(\alpha)=\e^{2\pi i\alpha}$ and $\e_k(\alpha)=\e(\alpha/k)$; $\bar x$ denotes an inverse of $x$ modulo the relevant modulus;
$\eps$ is an arbitrarily small positive constant, and $\delta\in(0,\frac1{20}]$ is fixed. Constants implied by $\ll$, $\asymp$
and $O$ may depend on $Q$, $\eps$ and $\delta$, and on further fixed data where indicated.
\section{Outline of the proof}\label{sec:strategy}

\begin{figure}[t]
\centering
\begin{tikzpicture}[>=stealth,node distance=6mm,
  box/.style={draw,rounded corners=2pt,align=flush center,inner sep=4pt,text width=8.6cm,font=\small},
  aux/.style={draw,dashed,rounded corners=2pt,align=left,inner sep=3pt,text width=3.9cm,font=\scriptsize},
  arr/.style={->,thick}]
\node[box,fill=green!8] (z) {\textbf{Zeros:} $Z$ has a zero $\tfrac12+i\gamma$\\[3pt] with \mbox{$|\gamma-T|\le H=T^{1+2\varepsilon}/K$}};
\node[box,below=of z,fill=blue!6] (p) {\textbf{Step 1} (\S\ref{subsec:zerodet}): critical Dirichlet polynomial\\[2pt]
\mbox{$\mathcal P(t)=\sum_n\eta(n)c(n)n^{-1/2-it}\ll(\log T)^{-2}$}\\[2pt]
for $|t-T|\le H$; $\eta$ of length $K$ around \mbox{$x_*=T\sqrt\Delta/2\pi$}};
\node[box,below=of p,fill=blue!6] (d) {\textbf{Step 2} (\S\S\ref{sec:FEV}--\ref{sec:transf}): twist by $e(hn/k)$, $h/k\approx1/\sqrt\Delta$;\\
Voronoi $+$ stationary phase $\Rightarrow$ sum over $w\in\mathbb Z^2$, $P=P(w)$:\\[2pt]
\mbox{$\mathcal P(t)\approx\sum_{\varsigma=\pm1}\sum_{P(w)\ll K}\beta(w)\,e\bigl(\varpi P+\varsigma\Phi(P)\bigr)\mathcal A_\varsigma(P)$}};
\node[box,below=of d,fill=blue!6] (s) {\textbf{Step 3} (\S\S\ref{sec:lemmas}--\ref{sec:expsum}):\\[2pt]
\mbox{$S(N)=\sum_{N\le P(w)\le N'}e(\cdots)\ll N^{2/3}K^{1/6}\log T$}\\[2pt]
$\Rightarrow$ admissible $K=T^{3/5-2\varepsilon}$, $H=T^{2/5+4\varepsilon}$};
\node[box,below=of s,fill=gray!8] (w) {two-dimensional Weyl differencing (Lemma~\ref{lem:weyl})\\ $+$ exact Hessian determinant (Lemma~\ref{lem:hess})\\
$+$ angular decomposition (Lemma~\ref{lem:decomp})\\ $+$ B-process and van der Corput (Lemma~\ref{lem:2d})};
\draw[arr] (p)--(z); \draw[arr] (d)--(p); \draw[arr] (s)--(d); \draw[arr] (w)--(s);
\node[aux,right=4mm of p] (a1) {Approximate formula, Ramachandra's lower bound, Stirling\\ \cite{JS,BSS}; Lemma~\ref{lem:approx}};
\node[aux,right=4mm of d] (a2) {Gauss sums (Lemma~\ref{lem:gauss}); functional equations and Voronoi formulas (Prop.~\ref{prop:FEdef}, Thm.~\ref{thm:FE}, Cor.~\ref{cor:unified}); stationary phase \cite{BKY} (Lemma~\ref{lem:statphase})};
\node[aux,right=4mm of s] (a3) {o-minimality (uniform finiteness); Graham--Kolesnik B-process};
\draw[->,gray] (a1)--(p); \draw[->,gray] (a2)--(d); \draw[->,gray] (a3)--(s);
\end{tikzpicture}
\caption{Structure of the proof of Theorems~\ref{thm:A}(a) and~\ref{thm:B}. Each box implies the box above it; dashed boxes
are the inputs.}
\label{fig:flow}
\end{figure}
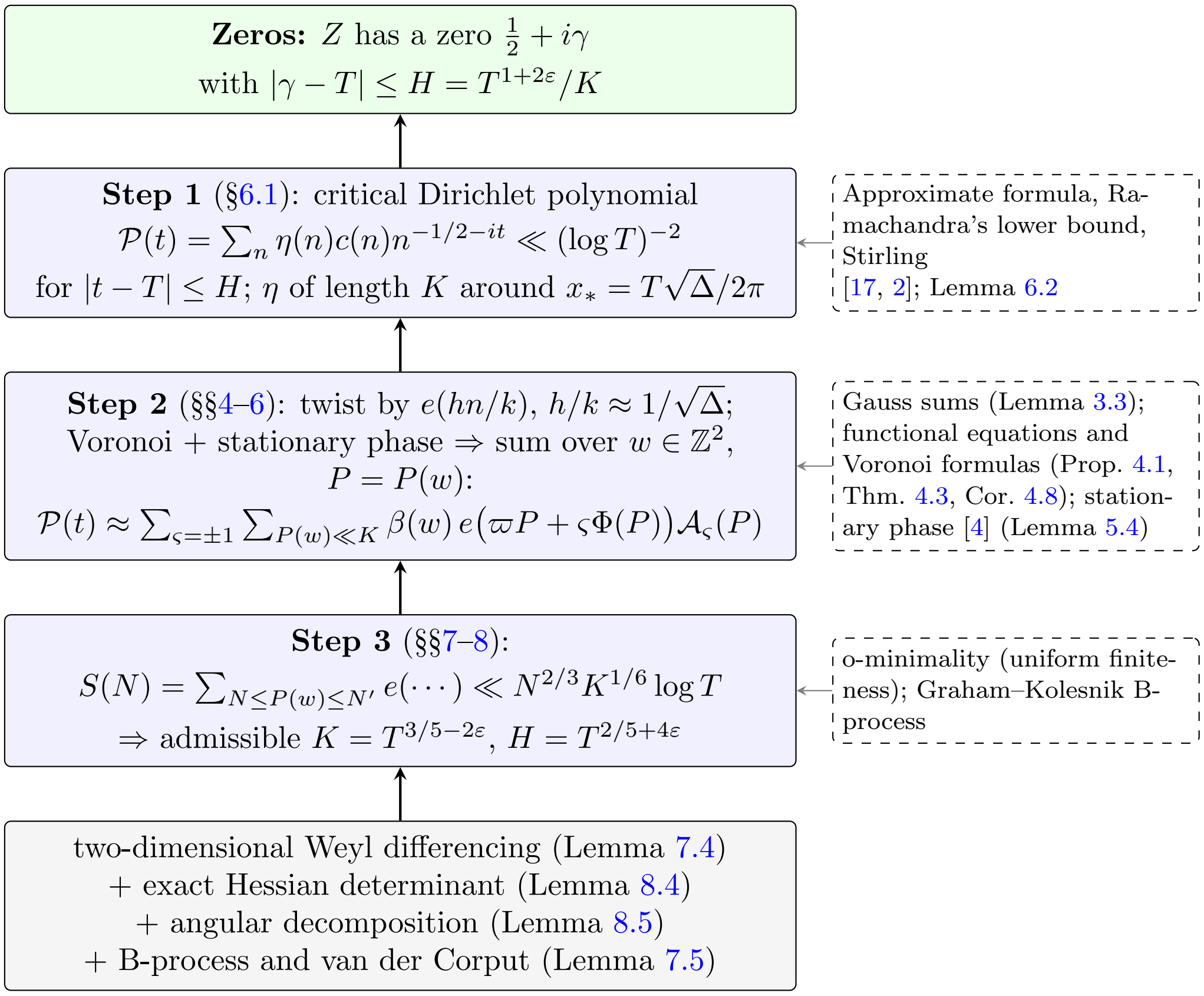

\subsection{The two cases in one language}\label{subsec:language}
Throughout, $Q=[a,b,c]$ is primitive (this is no loss of generality) and
\[
\Delta=|b^2-4ac|.
\]
We write $Z(s)=\sum_{n\ge1}c(n)n^{-s}$ for the zeta-function under consideration: $Z=\zeta_Q$ and $c=r_Q$ in the
definite case; $Z=\zQ$ or $Z=\Xi_Q^-$ and $c(n)=r_Q(n)+\omega r_Q(-n)$, $\omega=\pm1$, in the indefinite case.

\begin{figure}[t]
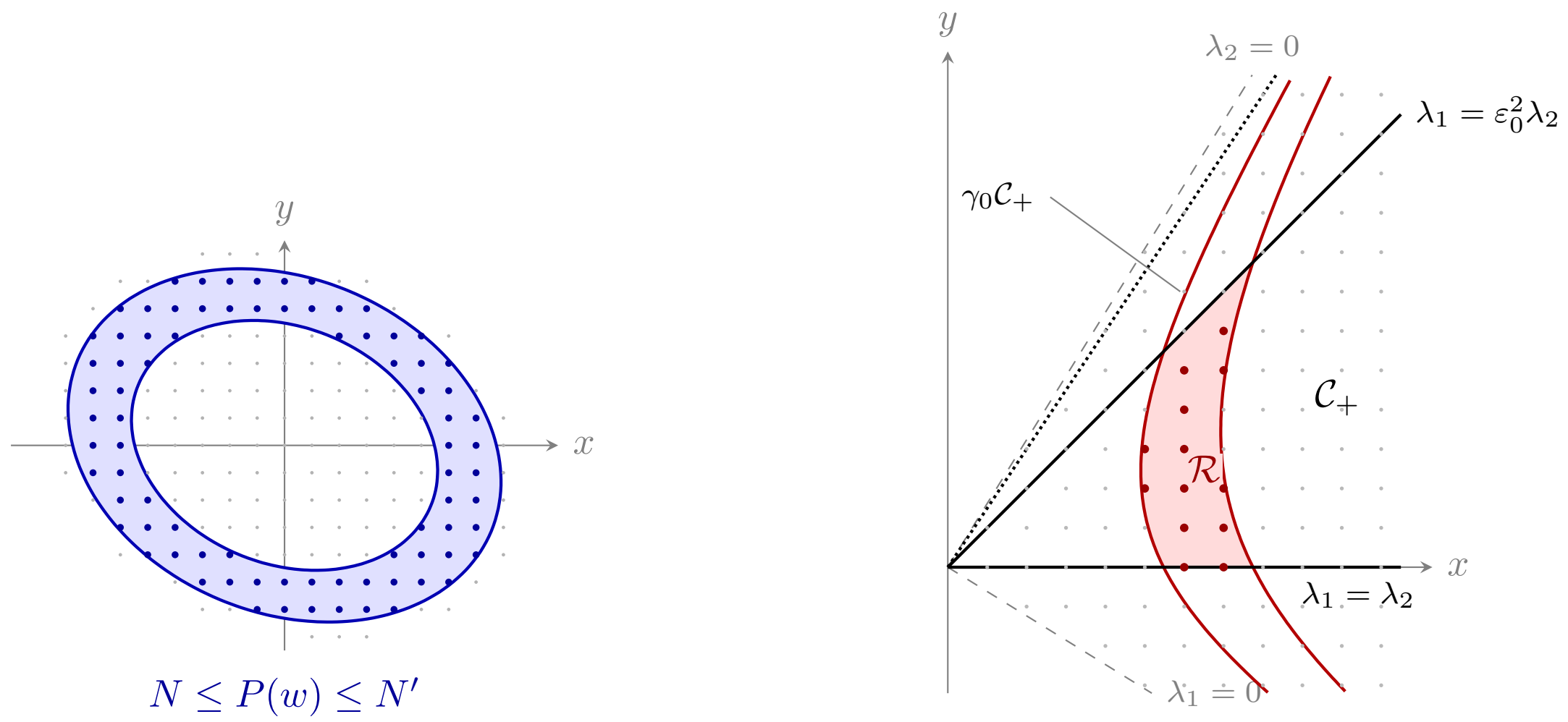

\centering
\begin{minipage}[b]{0.43\textwidth}\centering\input{figs/fig_annulus.tex}\end{minipage}\hfill
\begin{minipage}[b]{0.53\textwidth}\centering\input{figs/fig_sector.tex}\end{minipage}
\caption{The summation regions of $S(N)$. \emph{Left:} definite case, $P=2x^2+xy+3y^2$ ($d=-23$), the elliptic
annulus $60\le P(w)\le120$. \emph{Right:} indefinite case, $P=x^2+xy-y^2$ ($\Delta=5$, $\eps_0=\frac{3+\sqrt5}2$).
The fundamental cone $\mathcal C_+$ lies between the rays $\lambda_1=\lambda_2$ and $\lambda_1=\eps_0^2\lambda_2$; the dotted
ray bounds its image $\gamma_0\mathcal C_+$. The hyperbolic sector $\mathcal R$ is $\{w\in\mathcal C_+:30\le P(w)\le60\}$.
Large dots: the lattice points counted.}
\label{fig:regions}
\end{figure}

The central object after the transformation is the \emph{dual form}. Voronoi summation turns a sum over $n$ into a
sum over lattice points $w\in\Z^2$, and the frequency attached to $w$ is $P(w)$, where
\[
P=Q^*\ \text{(definite case)},\qquad P=|Q^*|=\pm Q^*\ \text{on the cone } \mathcal C^*_\pm\ \text{(indefinite case)},
\]
where $Q^*(w)=cw_1^2-bw_1w_2+aw_2^2$ is the adjoint form.
Here $\mathcal C^*_\pm$ are fundamental cones of the automorph group of $Q^*$ in $\{\pm Q^*>0\}$
(Lemma~\ref{lem:sector}). So in both cases $P$ is a positive quadratic form on the summation region, which is an
elliptic annulus in the definite case and a ``hyperbolic sector'' in the indefinite case (Figure~\ref{fig:regions}).
Table~\ref{tab:cases} collects the differences.

\begin{table}[t]
\centering\small
\renewcommand{\arraystretch}{1.25}
\begin{tabular}{@{}>{\raggedright}p{3.0cm}>{\raggedright}p{4.9cm}>{\raggedright\arraybackslash}p{5.6cm}@{}}
\toprule
& positive definite $Q$ & indefinite $Q$ \\
\midrule
$\Delta$ & $4ac-b^2$, not a square (Thm.~\ref{thm:A}(a)) & $b^2-4ac$, not a square\\
$r_Q(n)$ counts & solutions $v\in\Z^2$ of $Q(v)=n$ & $\Ga$-orbits of solutions, $n\in\Z\setminus\{0\}$\\
$Z(s)$, $c(n)$ & $\zeta_Q(s)$, $r_Q(n)$ & $\zQ$ or $\Xi^-_Q$, $r_Q(n)+\omega r_Q(-n)$\\
gamma factor & $(\sqrt\Delta/2\pi)^s\Ga(s)$ & $(\sqrt\Delta/\pi)^s\Ga(\frac{s+p}2)^2$, $p\in\{0,1\}$\\
Voronoi kernel & $J_0$ (Prop.~\ref{prop:Vdef}) & $-Y_0$ and $K_0$ (Thm.~\ref{thm:voronoi})\\
$\theta_\vs$ in Cor.~\ref{cor:unified} & $1$ & $i\vs$\\
dual set $\Dc$ & $\Z^2\setminus\{0\}$ & $(\mathcal C^*_+\cup\mathcal C^*_-)\cap\Z^2$ ($\Ga^*$-orbits)\\
dual form $P$ & $Q^*$ & $|Q^*|$\\
weights $\beta(w)$ & $\GQ(k,-h)/k$ & $\GQ(k,-h)/k$ if $Q^*(w)<0$, $\omega\GQ(k,h)/k$ if $Q^*(w)>0$\\
region of $S(N)$ & elliptic annulus & hyperbolic sector in a cone\\
\bottomrule
\end{tabular}
\caption{The two cases. In both, $|\beta(w)|=1$ when $(k,2\Delta)=1$.}
\label{tab:cases}
\end{table}

\subsection{Step 1: zero detection}\label{subsec:step1}
Let $W(t)$ be the Hardy function of $Z$: the function $Z(\frac12+it)$ multiplied by the gamma factor and an exponential
which make it real. If $W$ has no zero in $[T-H,T+H]$, then the integral
\[
I=\int_{-H}^HW(T+u)\e^{-(u/H_0)^2}du,\qquad H_0=HT^{-\eps},
\]
is large by Ramachandra's lower bound. On the other hand, an approximate formula for $Z(\frac12+it)$ and Stirling's
formula show that only the terms $n$ near the resonance point
\[
x_*=\frac{T\sqrt\Delta}{2\pi}
\]
contribute to $I$, since the other terms oscillate in $u$. This gives (Proposition~\ref{prop:zerodet}):
\begin{equation}\label{eq:strat1}
\text{if}\quad\Pc(t)=\sum_n\eta(n)c(n)n^{-1/2-it}\ll(\log T)^{-2}\quad(|t-T|\le H),\quad\text{then $Z$ has a zero in }[T-H,T+H].
\end{equation}
Here $\eta$ is a weight with $\eta=1$ on $[x_*-\frac K2,x_*+\frac K2]$, $\supp\eta\subset[x_*-K,x_*+K]$, and
\[
HK=T^{1+2\eps}.
\]
The trivial bound $\Pc(t)\ll KT^{-1/2}$ gives $K=T^{1/2-\eps}$, i.e.\ the exponent $\frac12$. Every power saving in
$\Pc$ beyond this lowers the exponent.

\subsection{Step 2: transformation of $\Pc(t)$}\label{subsec:step2}
\emph{Twist.} Choose $h/k$ close to $1/\sqrt\Delta$ with $(k,2\Delta)=1$ and $k\asymp\sqrt{T/K}$
(Lemma~\ref{lem:hk}), and write $n^{-it}=\e(-hn/k)\cdot\e(hn/k)n^{-it}$. The point of the twist is that the phase
\[
g(x)=\frac hkx-\frac t{2\pi}\log x
\]
of $\e(hx/k)x^{-it}$ is stationary at $x_0=tk/(2\pi h)$, which lies in the middle half of $\supp\eta$. Then
\[
\Pc(t)=\sum_nc(n)\e\Bigl(-\frac{hn}k\Bigr)f(n),\qquad f(x)=\eta(x)x^{-1/2}\e(g(x)).
\]
\emph{Voronoi.} By the twisted Voronoi formulas (Corollary~\ref{cor:unified}), with $q=k\sqrt\Delta$,
\[
\Pc(t)=\mathfrak m\int f+\frac\pi q\sum_{w\in\Dc}\beta(w)\,\e_k\bigl(\hb\Db P(w)\bigr)\sum_{\vs=\pm1}\theta_\vs
\int_0^\infty f(x)H_0^{(\vs)}\Bigl(\frac{4\pi\sqrt{P(w)x}}q\Bigr)dx+(\text{$K_0$-terms}),
\]
where $H_0^{(+1)}=H^{(1)}_0$ and $H_0^{(-1)}=H^{(2)}_0$ are the Hankel functions, and $\mathfrak m\ll1/k$.

\emph{Stationary phase.} With $H_0^{(\vs)}(z)\approx\sqrt{2/\pi z}\,\e^{i\vs(z-\pi/4)}$, the $w$-th integral has the phase
$g(x)+2\vs\sqrt{P(w)x}/q$. Its stationary point $x_\vs(n)$, $n=P(w)$, moves away from $x_0$ as $n$ grows, and lies in
$\supp\eta$ only if $n\ll N_0\asymp K$ (Figure~\ref{fig:statpoints}). At the stationary point the phase equals, up to a
constant independent of $n$,
\[
-\frac n{2q^2r}+\vs\Phi(n),\qquad r=\frac hk,\qquad \Phi(u)=\frac t\pi\,\phi\Bigl(\frac{\pi u}{2q^2rt}\Bigr)\approx2\kappa\sqrt u,\qquad\kappa\asymp K^{1/2},
\]
with $\phi(x)=\arsinh(x^{1/2})+(x+x^2)^{1/2}$ (Lemma~\ref{lem:statphase}). The rational phase $\hb\Db P(w)/k$ and the term
$-P(w)/(2q^2r)$ combine into $\varpi P(w)$, with $\varpi=\hb\Db/k-1/(2q^2r)$. The result is (Proposition~\ref{prop:transf})
\begin{equation}\label{eq:strat2}
\Pc(t)=\e(\cdots)\sum_{\vs=\pm1}\sum_{w\in\Dc}\beta(w)\,\e\bigl(\varpi P(w)+\vs\Phi(P(w))\bigr)\mathcal A_\vs(P(w))
+O\bigl(T^{1/2}K^{-1}+KT^{-5/6}+(K/T)^{1/2}\bigr),
\end{equation}
where the amplitudes vanish for $u>5N_0$ and satisfy
\[
\mathcal A_\vs(u)\ll K^{1/4}u^{-1/4}T^{-1/2},\qquad\Var_{[N,2N]}\mathcal A_\vs\ll K^{1/4}N^{-1/4}T^{-1/2}.
\]
The error terms come from the lower order terms of the stationary phase expansion, which are smaller by
$T/K^2+T^{-1/3}$ and are estimated trivially. They are $\ll T^{-\delta}$ in the range \eqref{eq:Krange} below.

\emph{Partial summation} in $u=P(w)$ over dyadic ranges now reduces \eqref{eq:strat1} to the following bound
(Corollary~\ref{cor:reduction}): for $1\le N\le5N_0$, $N<N'\le2N$ and $\vs=\pm1$,
\begin{equation}\label{eq:strat3}
K^{1/4}N^{-1/4}T^{-1/2}\,|S(N)|\ll T^{-\eps},\qquad S(N)=\sum_{\substack{w\in\Dc_j\\N\le P(w)\le N'}}\e\bigl(\varpi P(w)+\vs\Phi(P(w))\bigr),
\end{equation}
where $\Dc_j$ is $\Dc$ (definite case) or one of the two cones (indefinite case), on which $\beta$ is constant.

\subsection{Step 3: the double exponential sum}\label{subsec:step3}
The sum $S(N)$ has $\asymp N$ terms, $|w|\asymp L=N^{1/2}$, and phase of size $\kappa L$. The trivial bound $|S|\ll N$ in
\eqref{eq:strat3} gives $K=T^{1/2}$ again. We prove (Theorem~\ref{thm:SN})
\begin{equation}\label{eq:strat4}
S(N)\ll N^{2/3}K^{1/6}\log T\qquad(c_0K^{1/2}\le N\ll K,\ T^{1/2+\delta}\le K\le T^{3/5}).
\end{equation}
Inserting this into \eqref{eq:strat3} gives $K^{5/12}N^{5/12}T^{-1/2}\log T\le K^{5/6}T^{-1/2}\log T$, which is small for
$K=T^{3/5-2\eps}$. Hence $H=T^{2/5+4\eps}$. The proof of \eqref{eq:strat4} has four steps.

\begin{figure}[t]
\centering
\input{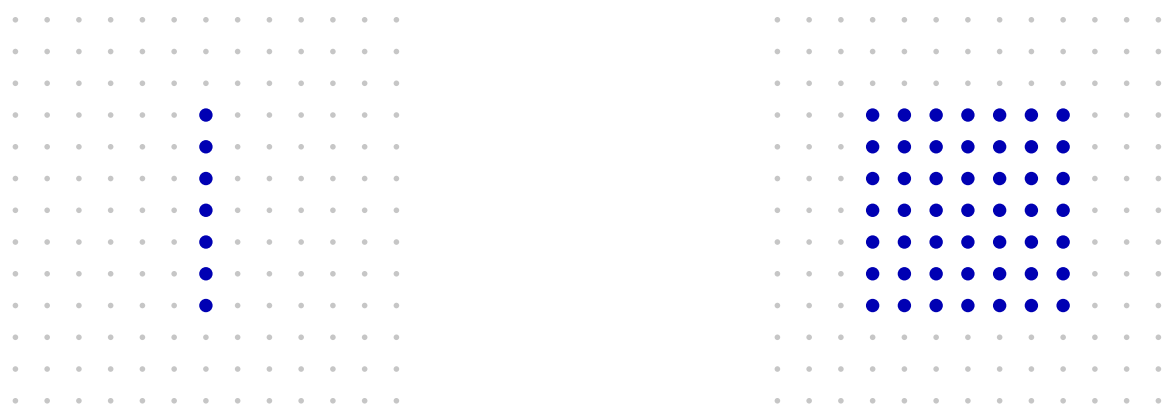}
\caption{Shift sets in Weyl's inequality. \emph{Left:} one-dimensional differencing as in \cite{BSS}, with $\asymp M$
shifts and diagonal loss $N^2/M$. \emph{Right:} the two-dimensional differencing used here, with $\asymp M^2$ shifts and
diagonal loss $N^2/M^2$ (here $M=4$).}
\label{fig:shifts}
\end{figure}

\emph{(i) Weyl differencing in two variables} (Lemma~\ref{lem:weyl}). With shifts $\ell\in\Z^2$, $|\ell|_\infty<M$
(Figure~\ref{fig:shifts}),
\[
|S(N)|^2\ll\frac{L^2}{M^2}\Bigl(L^2+\sum_{0<|\ell|\ll M}|S_\ell|\Bigr),\qquad
S_\ell=\sum_w\e\bigl(f(w+\ell)-f(w)\bigr).
\]
The quadratic part $\varpi P$ of the phase becomes linear in $w$, so the rational number $\hb\Db/k$ disappears from all
second derivatives. The relevant part of the differenced phase is
\[
g_\ell(w)=\Phi(P(w+\ell))-\Phi(P(w))\approx2\kappa\bigl(\sqrt{P(w+\ell)}-\sqrt{P(w)}\bigr).
\]
\emph{(ii) The Hessian} (Lemma~\ref{lem:hess}). The Hessian of $\sqrt P$ has rank one,
\[
\nabla^2\sqrt P(w)=-\frac d4P(w)^{-3/2}(Jw)(Jw)^T,\qquad J=\begin{pmatrix}0&-1\\1&0\end{pmatrix},\quad d=\operatorname{disc}P,
\]
and therefore the Hessian $\Hs_\ell$ of $g_\ell$ satisfies, with $\psi$ the angle between $w$ and $\ell$,
\[
\|\Hs_\ell\|\asymp\lambda:=\frac{\kappa|\ell|}{L^2},\qquad|\det\Hs_\ell|\asymp(\lambda\sin\psi)^2 .
\]
The determinant degenerates only when $w$ is nearly parallel to $\ell$.

\begin{figure}[t]
\centering
\begin{tikzpicture}[>=stealth,scale=0.62]
\fill[blue!8,even odd rule] (0,0) circle (5.6) (0,0) circle (4.0);
\draw[blue!60!black] (0,0) circle (4.0) (0,0) circle (5.6);
\fill[gray!55] (25.0:4.0) arc (25.0:31.0:4.0) -- (31.0:5.6) arc (31.0:25.0:5.6) -- cycle;
\fill[orange!45] (37.0:4.0) arc (37.0:46.0:4.0) -- (46.0:5.6) arc (46.0:37.0:5.6) -- cycle;
\fill[orange!45] (10.0:4.0) arc (10.0:19.0:4.0) -- (19.0:5.6) arc (19.0:10.0:5.6) -- cycle;
\fill[gray!55] (205.0:4.0) arc (205.0:211.0:4.0) -- (211.0:5.6) arc (211.0:205.0:5.6) -- cycle;
\fill[orange!45] (217.0:4.0) arc (217.0:226.0:4.0) -- (226.0:5.6) arc (226.0:217.0:5.6) -- cycle;
\fill[orange!45] (190.0:4.0) arc (190.0:199.0:4.0) -- (199.0:5.6) arc (199.0:190.0:5.6) -- cycle;
\draw[dashed] (208.0:6.3)--(28.0:6.3) node[above right]{\small $\mathbb R\ell$};
\draw[->,very thick] (0,0)--(28.0:1.0) node[above left=-1pt]{\small $\ell$};
\fill (0,0) circle (1.2pt) node[below left]{\small $0$};
\draw[green!45!black] (3.658,2.080)--(2.851,3.699);
\draw[green!45!black,line width=1.6pt] (3.413,2.572)--(3.096,3.206);
\draw[green!45!black] (4.051,2.375)--(3.207,4.067);
\draw[green!45!black,line width=1.6pt] (3.805,2.867)--(3.452,3.575);
\draw[green!45!black] (4.404,2.641)--(3.527,4.399);
\draw[green!45!black,line width=1.6pt] (4.158,3.133)--(3.773,3.907);
\draw[->,thick,green!45!black] (4.6,-5.2)--++(-0.491,0.984) node[right]{\small $e$};
\node at (4.981,4.406) {\small $\mathcal A(\Psi)$};
\draw[gray!70!black] (5.100,2.655)--(6.404,1.597) node[right]{\small $\psi\le\Psi_0$};
\node[blue!60!black] at (0,-4.8) {\small $|w|\asymp L$};
\draw (28.0:2.3) arc (28.0:37.0:2.3); \draw (28.0:2.9) arc (28.0:46.0:2.9);
\node at (32.05:2.0) {\scriptsize $\Psi$}; \node at (49.0:3.2) {\scriptsize $2\Psi$};
\end{tikzpicture}
\caption{The decomposition of Lemma~\ref{lem:decomp} (definite case, schematic). The double cone $\psi\le\Psi_0$ around the
line $\R\ell$ (grey) is counted trivially. The four pieces of $\mathcal A(\Psi)$, where the angle between $w$ and $\R\ell$
lies in $(\Psi,2\Psi]$, are shaded. Lines parallel to a direction $e$ transverse to $w$ meet such a piece in segments of
length $\ll L\Psi$ (bold); on them the B-process is applied, and van der Corput's estimate in the transverse variable.}
\label{fig:decomp}
\end{figure}
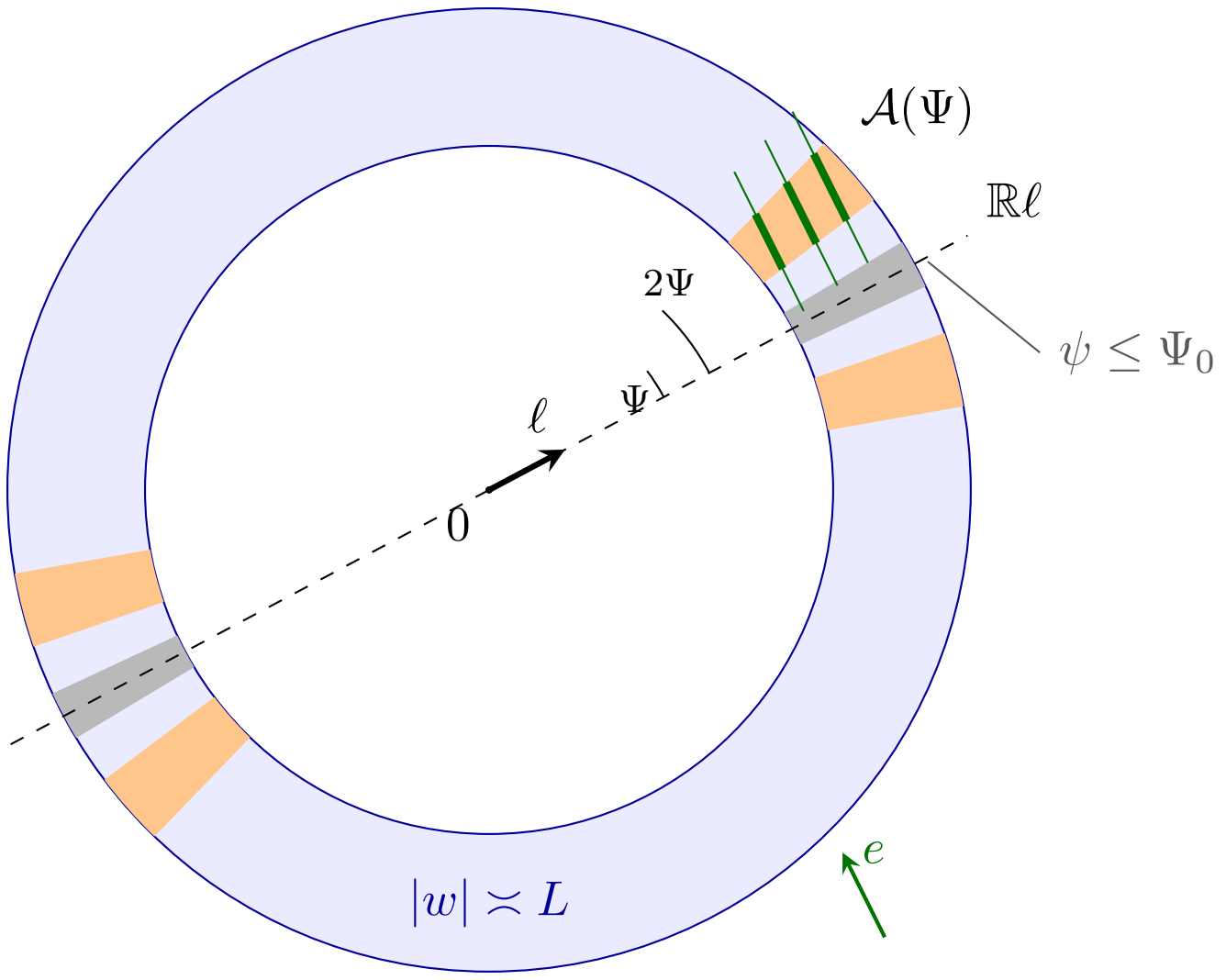

\emph{(iii) Decomposition and summation} (Lemmas~\ref{lem:decomp}, \ref{lem:2d}). The points with $\psi\le\Psi_0=(\kappa|\ell|)^{-1/2}$
are counted trivially. The rest is split dyadically into pieces $\psi\asymp\Psi$ (Figure~\ref{fig:decomp}). On each piece
we apply the B-process in a direction $e$ transverse to $w$, and then van der Corput's second derivative estimate in the
other variable, where the second derivative is $\det\Hs_\ell/\partial_e^2g_\ell\asymp\lambda\Psi^2$. The result is
\[
|S_\ell|\ll\Bigl(\kappa|\ell|+L+\frac{L^2}{(\kappa|\ell|)^{1/2}}\Bigr)\log^2T .
\]
\emph{(iv) Choice of $M$.} Summing over $\ell$,
\[
|S(N)|^2\ll\Bigl(\frac{L^4}{M^2}+L^2\kappa M+L^3+\frac{L^4}{(\kappa M)^{1/2}}\Bigr)\log^2T,
\]
and $M\asymp(L^2/\kappa)^{1/3}$ balances the first two terms. This gives $|S(N)|\ll L^{4/3}\kappa^{1/3}\log T$, i.e.\ \eqref{eq:strat4}.

\subsection{Bookkeeping of exponents}
At $N\asymp K$ the sum $S(N)$ has $\asymp K$ terms and phase of size $\asymp K$, and a bound $|S(K)|\ll K^{1-\sigma}$ in
\eqref{eq:strat3} allows $K=T^{1/(2(1-\sigma))}$, i.e.\ the exponent \eqref{eq:theta}. Theorem~\ref{thm:SN} gives $\sigma=\frac16$,
hence $K=T^{3/5-2\eps}$ and $\theta=\frac25$.

\subsection{Square $\Delta$}\label{subsec:squareintro}
If $\Delta=f^2$ (then $f=2m$ is even), the critical ratio $1/\sqrt\Delta=1/f$ is rational, and we may take $h/k=1/f$ itself.
Then $q=f^2$ is bounded, and the dual sum is short: $N_0\asymp K^2/T$. As for Jutila's cusp forms \cite{J1}, this gives the
exponent $1/3$, provided the main term $\mathfrak m\int f$, which is $\asymp\GQ(f,-1)$, vanishes. This Gauss sum vanishes if
and only if $m$ is odd (Lemma~\ref{lem:Gsq}). For even $m$ the main term contributes a non-zero constant to the Hardy
function itself (Proposition~\ref{prop:square}), so the method cannot work; for $x^2+4y^2$ the zeros are nevertheless
spaced $O(1)$ apart (Proposition~\ref{prop:x24y2}).

\begin{table}[t]
\centering\small
\renewcommand{\arraystretch}{1.2}
\begin{tabular}{@{}>{$}l<{$}>{\raggedright}p{6.6cm}>{$}l<{$}l@{}}
\toprule
\text{symbol} & meaning & \text{size} & defined\\
\midrule
T & height & & \\
K & length of $\supp\eta$ & T^{3/5-2\eps} & \eqref{eq:Krange}\\
H & half-length of the interval for $\gamma$; $HK=T^{1+2\eps}$ & T^{2/5+4\eps} & \S\ref{subsec:zerodet}\\
t,\ \tau & $|t-T|\le H$; $\tau=t/2\pi$ & T & \\
x_* & resonance point $T\sqrt\Delta/2\pi$ & T & \eqref{eq:xstar}\\
h/k=r & approximation to $1/\sqrt\Delta$, $(k,2\Delta)=1$ & k\asymp\sqrt{T/K} & Lemma~\ref{lem:hk}\\
q & $k\sqrt\Delta$, the scale in the Bessel functions & \sqrt{T/K} & \eqref{eq:qdef}\\
x_0 & $\tau/r$, stationary point of $g$; $|x_0-x_*|\le K/2$ & T & \eqref{eq:x0}\\
N_0 & $r^2q^2K^2/x_*$, effective length of the dual sum & K & Lemma~\ref{lem:statphase}\\
\varpi & $\hb\Db/k-1/(2q^2r)$ & & \eqref{eq:varpi}\\
\kappa & $(t/(2\pi q^2r))^{1/2}$, so $\Phi(u)\approx2\kappa\sqrt u$ & K^{1/2} & \eqref{eq:Phi}\\
N,\ L & dyadic range of $P(w)$; $L=N^{1/2}$ & N\le C_1K & \S\ref{sec:expsum}\\
M & number of shifts per coordinate & (L^2/\kappa)^{1/3} & \S\ref{subsec:proofSN}\\
\ell,\ \lambda & shift; $\lambda=\kappa|\ell|L^{-2}\asymp\|\Hs_\ell\|$ & |\ell|\ll M & \S\ref{subsec:diffphase}\\
\psi,\ \Psi_0 & angle between $w$ and $\ell$; $\Psi_0=(\kappa|\ell|)^{-1/2}$ & & \S\ref{subsec:decomp}\\
\bottomrule
\end{tabular}
\caption{Parameters, with their sizes in the proofs of Theorems~\ref{thm:A}(a) and~\ref{thm:B}.}
\label{tab:parameters}
\end{table}

\section{Binary forms: adjoint form, automorphs, Gauss sums}\label{sec:forms}

This section collects the algebraic facts used later: the adjoint form, which carries the frequencies of the dual
sums; for indefinite forms, the group of automorphs and its fundamental cones, which replace the finite sums of the
definite case by sums over cones; and the quadratic Gauss sums, which appear as coefficients in the Voronoi formulas.

\subsection{Matrices and the adjoint form}
For $Q=[a,b,c]$ put
\[
A=\begin{pmatrix}2a&b\\ b&2c\end{pmatrix},\qquad A^\dagger=\begin{pmatrix}2c&-b\\ -b&2a\end{pmatrix},\qquad
Q(v)=\tfrac12v^TAv,\qquad Q^*(w)=\tfrac12w^TA^\dagger w .
\]
Then $AA^\dagger=(\det A)I$ with $\det A=4ac-b^2$, i.e.\ $\det A=\Delta$ for definite and $\det A=-\Delta$ for indefinite $Q$.
The bilinear form is $B(v,w)=v^TAw=Q(v+w)-Q(v)-Q(w)$. The adjoint form $Q^*=[c,-b,a]$ has the same discriminant as $Q$,
and $Q^*(w_2,-w_1)=Q(w_1,w_2)$, so $Q^*$ is properly equivalent to $Q$.

\subsection{Indefinite forms: automorphs and the fundamental cone}\label{subsec:automorphs}
Let $Q$ be indefinite with discriminant $\Delta>0$ not a square. Factor $Q=\lambda_1\lambda_2$ over $\R$ with real linear
forms, and let $\Theta$ be the matrix with rows $\lambda_1,\lambda_2$. Then
\[
A=2\Theta^T\begin{pmatrix}0&1/2\\1/2&0\end{pmatrix}\Theta,\qquad(\det\Theta)^2=\Delta .
\]
Let $t_0^2-\Delta u_0^2=4$ be the fundamental solution and $\eps_0=(t_0+u_0\sqrt\Delta)/2$. Then $\Ga=\{\pm\gamma_0^j:j\in\Z\}$, where
\[
\lambda_1\circ\gamma_0=\eps_0\lambda_1,\qquad\lambda_2\circ\gamma_0=\eps_0^{-1}\lambda_2 .
\]
No $\gamma\ne1$ in $\Ga$ fixes a non-zero vector. The dual linear forms $(\mu_1,\mu_2)^T=\Theta^{-T}w$ satisfy
\begin{equation}\label{eq:mumu}
\mu_1(w)\mu_2(w)=\tfrac12w^TA^{-1}w=-\frac{Q^*(w)}\Delta .
\end{equation}
With $\gamma^*=\gamma^{-T}$, the group $\Ga^*=\{\gamma^*:\gamma\in\Ga\}=SO_{Q^*}(\Z)$ acts by $\mu_1\circ\gamma^*=\eps^{-1}\mu_1$ and
$\mu_2\circ\gamma^*=\eps\mu_2$ whenever $\lambda_1\circ\gamma=\eps\lambda_1$.

\begin{lemma}[Fundamental cone; see Figure~\ref{fig:regions}]\label{lem:sector}
The set
\[
\mathcal C_+=\Bigl\{v:\lambda_1(v)>0,\ \lambda_2(v)>0,\ 1\le\frac{\lambda_1(v)}{\lambda_2(v)}<\eps_0^2\Bigr\}
\]
is a fundamental domain for $\Ga$ on $\{Q>0\}$, and
\[
\mathcal C_-=\Bigl\{v:\lambda_1(v)>0>\lambda_2(v),\ 1\le\Bigl|\frac{\lambda_1(v)}{\lambda_2(v)}\Bigr|<\eps_0^2\Bigr\}
\]
is one for $\{Q<0\}$. Their closures minus the origin lie in the open cones $\{\pm Q>0\}$, so $|v|\asymp|Q(v)|^{1/2}$ on
$\mathcal C_\pm$. In particular $r_Q(\pm n)=\#\{v\in\mathcal C_\pm\cap\Z^2:Q(v)=\pm n\}$, and
$\#\{v\in\mathcal C_\pm\cap\Z^2:|Q(v)|\le X\}\ll X$ for $X\ge1$.
\end{lemma}
\begin{proof}
The element $-1$ interchanges the two components of $\{Q>0\}$, and $\gamma_0$ multiplies $\lambda_1/\lambda_2$ by $\eps_0^2$.
The ratio $\lambda_1/\lambda_2$ is bounded above and below on $\mathcal C_+$, and similarly on $\mathcal C_-$. The last claim
holds because $\mathcal C_\pm\cap\{|Q|\le X\}$ is contained in a disc of radius $\ll X^{1/2}$.
\end{proof}

Applied to $Q^*$, Lemma~\ref{lem:sector} gives fundamental cones $\mathcal C^*_\pm$ for $\Ga^*$ on $\{\pm Q^*>0\}$.

\begin{remark}[Ideal classes]\label{rem:ideals}
Let $\Delta$ be fundamental and $\mathfrak a=[a,\frac{b-\sqrt\Delta}2]$, so that
$Q(x,y)=N(xa+y\frac{b-\sqrt\Delta}2)/N\mathfrak a$ (cf.\ \cite[Prop.~6]{BR}). Then $\Ga\backslash\{Q>0\}$ corresponds to
$\{\alpha\in\mathfrak a:\alpha\gg0\}/U^+$. Thus $Z_Q(s)=\sum_{n\ge1}r_Q(n)n^{-s}$ is the zeta-function of the narrow class of
$\mathfrak a^{-1}$. Correspondingly, $\zQ$ is the zeta-function of the wide class if $N(U)=\{1\}$, and twice it otherwise;
in the latter case $Q\sim-Q$.
\end{remark}

\subsection{Gauss sums}\label{subsec:gauss}
In this subsection $Q=[a,b,c]$ is any integral binary form with $\det A\ne0$. For $k\ge1$, $(h,k)=1$ and $w\in\Z^2$ put
\begin{gather*}
\GQ(k,h,w)=\sum_{x\in(\Z/k\Z)^2}\e_k\bigl(hQ(x)+x\cdot w\bigr),\qquad\GQ(k,h)=\GQ(k,h,0),\\
\mathcal R_k=\{c\bmod k:Ac\equiv0\pmod k\}.
\end{gather*}

\begin{lemma}\label{lem:gauss}
\begin{enumerate}
\item[(a)] If $Q$ is indefinite, $\GQ(k,h,\gamma^*w)=\GQ(k,h,w)$ for $\gamma\in\Ga$. In general
$\GQ(k,h,-w)=\GQ(k,h,w)$ and $\GQ(k,-h,w)=\overline{\GQ(k,h,w)}$.
\item[(b)] $\GQ(k,h,w)=0$ unless $hQ(c)+c\cdot w\equiv0\pmod k$ for all $c\in\mathcal R_k$, in which case
$|\GQ(k,h,w)|=k|\mathcal R_k|^{1/2}$. If $Q$ is primitive, then $|\mathcal R_k|\le2(k,\det A)$.
\item[(c)] If $(k,2\det A)=1$, then $|\GQ(k,h)|=k$ and
\[
\GQ(k,h,w)=\GQ(k,h)\,\e_k\bigl(-\hb\,\overline{\det A}\,Q^*(w)\bigr).
\]
\end{enumerate}
\end{lemma}
\begin{proof}
(a) Substitute $x=\gamma y$, respectively $x=-y$; the last identity is clear.

(b) Expanding $|\GQ|^2$ and putting $x=y+c$ gives
\[
|\GQ|^2=\sum_c\e_k\bigl(hQ(c)+c\cdot w\bigr)\sum_y\e_k\bigl(hy^TAc\bigr)=k^2\sum_{c\in\mathcal R_k}\e_k\bigl(hQ(c)+c\cdot w\bigr),
\]
using $(h,k)=1$. On $\mathcal R_k$ the map $c\mapsto hQ(c)+c\cdot w$ is additive modulo $k$, since $Q(c+c')=Q(c)+Q(c')+c^TAc'$.
So the last sum is $|\mathcal R_k|$ or $0$. If $\mathrm{diag}(d_1,d_2)$ is the Smith form of $A$, then $d_1d_2=|\det A|$,
$d_1=\gcd(2a,b,2c)$ divides $2$ for primitive $Q$, and $|\mathcal R_k|=(d_1,k)(d_2,k)\le2(k,\det A)$.

(c) Here $\mathcal R_k=\{0\}$, so $|\GQ(k,h)|=k$ by (b). Take $c\equiv-\hb\,\overline{\det A}\,A^\dagger w\pmod k$. Then
$hAc\equiv-w$, hence
\[
hQ(x-c)=hQ(x)-hx^TAc+hQ(c)\equiv hQ(x)+x\cdot w+hQ(c)\pmod k .
\]
Moreover $2Q(c)=c^TAc\equiv-\hb\,c\cdot w=\hb^2\,\overline{\det A}\,w^TA^\dagger w=2\hb^2\,\overline{\det A}\,Q^*(w)$, and $2$ is
invertible modulo $k$. Summing over $x$ gives $\GQ(k,h,w)=\GQ(k,h)\e_k(-hQ(c))=\GQ(k,h)\e_k(-\hb\,\overline{\det A}\,Q^*(w))$.
\end{proof}

\section{Functional equations and Voronoi formulas}\label{sec:FEV}

For $\sigma>1$ we put
\[
L_Q(s,h/k)=\sum_{n\ge1}r_Q(n)\e_k(hn)n^{-s}.
\]
In both cases the functional equation comes from Poisson summation over the shifted lattices $a+k\Z^2$, and the Voronoi
formula follows by Mellin inversion.

\subsection{Positive definite forms}

\begin{proposition}\label{prop:FEdef}
Let $Q$ be positive definite, $k\ge1$ and $(h,k)=1$. Then $L_Q(s,h/k)$ continues meromorphically to $\C$. Its only
possible singularity is a simple pole at $s=1$, with residue $2\pi\GQ(k,h)/(k^2\sqrt\Delta)$. Moreover
\begin{equation}\label{eq:FEdef}
(2\pi)^{-s}\Ga(s)L_Q(s,h/k)=\Delta^{\frac12-s}k^{-2s}(2\pi)^{s-1}\Ga(1-s)\sum_{w\in\Z^2\setminus\{0\}}\frac{\GQ(k,h,w)}{Q^*(w)^{1-s}}\qquad(\sigma<0).
\end{equation}
\end{proposition}
\begin{proof}
For $a\in\Z^2$ and $y>0$ let $\theta_a(y)=\sum_{v\in a+k\Z^2}\e^{-2\pi yQ(v)}$. Poisson summation in the form
\[
\sum_{v\in a+k\Z^2}F(v)=k^{-2}\sum_{w\in\Z^2}\e_k(a\cdot w)\widehat F(w/k)
\]
applies to $F(v)=\e^{-\pi y\,v^TAv}$, whose Fourier transform is $\widehat F(\xi)=(y\sqrt\Delta)^{-1}\e^{-\pi\xi^TA^{-1}\xi/y}$. Since
$A^{-1}=A^\dagger/\Delta$, we have $\xi^TA^{-1}\xi=2Q^*(\xi)/\Delta$, and therefore
\begin{equation}\label{eq:thetadef}
\theta_a(y)=\frac1{k^2\sqrt\Delta\,y}\sum_{w\in\Z^2}\e_k(a\cdot w)\,\e^{-2\pi Q^*(w)/(\Delta k^2y)} .
\end{equation}
Let $Z(s;a)=\sum_{0\ne v\in a+k\Z^2}Q(v)^{-s}$ and $Z^*(s;a)=\sum_{w\ne0}\e_k(a\cdot w)Q^*(w)^{-s}$. Then
\[
(2\pi)^{-s}\Ga(s)Z(s;a)=\int_0^\infty\bigl(\theta_a(y)-\mathbf 1_{a\equiv0}\bigr)y^{s}\frac{dy}y .
\]
Splitting the integral at $y=1/(k\sqrt\Delta)$ and applying \eqref{eq:thetadef} on the lower part, as in Riemann's proof,
gives the meromorphic continuation and
\[
(2\pi)^{-s}\Ga(s)Z(s;a)=\Delta^{\frac12-s}k^{-2s}(2\pi)^{s-1}\Ga(1-s)Z^*(1-s;a).
\]
The term $w=0$ in \eqref{eq:thetadef} produces the only pole of $Z(s;a)$, at $s=1$ with residue $2\pi/(k^2\sqrt\Delta)$; the
term $\mathbf 1_{a\equiv0}$ produces a pole at $s=0$, which is cancelled by $\Ga(s)$. Since $Q(v)\equiv Q(a)\pmod k$ for
$v\in a+k\Z^2$,
\[
L_Q(s,h/k)=\sum_{a\bmod k}\e_k(hQ(a))Z(s;a),
\]
and $\sum_{a\bmod k}\e_k(hQ(a))\e_k(a\cdot w)=\GQ(k,h,w)$. This proves the proposition.
\end{proof}

\begin{proposition}[Voronoi formula, definite case]\label{prop:Vdef}
Let $Q$ be positive definite, $k\ge1$, $(h,k)=1$, and let $f\in C^3$ have compact support in $(0,\infty)$. Then
\begin{multline*}
\sum_{n\ge1}r_Q(n)\e_k(hn)f(n)=\frac{2\pi\GQ(k,h)}{k^2\sqrt\Delta}\int_0^\infty f(x)\,dx\\
+\frac{2\pi}{k^2\sqrt\Delta}\sum_{w\ne0}\GQ(k,h,w)\int_0^\infty f(x)J_0\Bigl(\frac{4\pi\sqrt{Q^*(w)x}}{k\sqrt\Delta}\Bigr)dx .
\end{multline*}
\end{proposition}
\begin{proof}
Write the left side as $\frac1{2\pi i}\int_{(2)}L_Q(s,h/k)\tilde f(s)\,ds$, where $\tilde f(s)=\int f(x)x^{s-1}dx\ll(1+|s|)^{-3}$
in vertical strips. Shift the contour to $\sigma=-\frac18$, picking up the residue at $s=1$, which gives the first term. On
$\sigma=-\frac18$, \eqref{eq:FEdef} and Stirling's formula give $L_Q\ll(1+|s|)^{5/4}$, and $L_Q$ has finite order in vertical
strips, so the shift is justified by Phragm\'en--Lindel\"of. Insert \eqref{eq:FEdef}, written as
\[
L_Q(s,h/k)=\sum_{w\ne0}\GQ(k,h,w)\,\frac{\sqrt\Delta}{2\pi Q^*(w)}\,\frac{\Ga(1-s)}{\Ga(s)}\,y_w^{\,s},\qquad y_w=\frac{4\pi^2Q^*(w)}{\Delta k^2},
\]
and integrate term by term. For each $w$, move the line of integration to $\sigma=\frac12$ (there are no poles in
$-\frac18\le\sigma\le\frac12$), and use Mellin--Parseval together with
\[
\int_0^\infty J_0(2\sqrt u)\,u^{z-1}du=\frac{\Ga(z)}{\Ga(1-z)}\qquad(0<\Re z<\tfrac34).
\]
This gives $\frac1{2\pi i}\int\frac{\Ga(1-s)}{\Ga(s)}y^s\tilde f(s)\,ds=y\int_0^\infty f(x)J_0(2\sqrt{yx})\,dx$, and
$\frac{\sqrt\Delta}{2\pi Q^*(w)}\,y_w=\frac{2\pi}{k^2\sqrt\Delta}$. Absolute convergence of the resulting series follows from
$\int fJ_0(c\sqrt x)\,dx\ll_fc^{-7/2}$ (three integrations by parts, using $\frac d{dx}(\frac{2\sqrt x}cJ_1(c\sqrt x))=J_0(c\sqrt x)$
and its analogues) and $|\GQ(k,h,w)|\ll_k1$.
\end{proof}

\subsection{Indefinite forms}
Let $Q$ be indefinite, $k\ge1$ and $(h,k)=1$, and let $\rho(n)=\rho_Q(n;h/k)$ be as in \S\ref{sec:intro}. In this subsection
we prove Theorem~\ref{thm:voronoi}.

The dual coefficients $\rho(n)$ are well defined by
Lemma~\ref{lem:gauss}(a), and $|\rho(n)|\le k\sqrt{2(k,\Delta)}\,r_{Q^*}(n)$. Put
\[
D_\pm(s)=\sum_{n\ge1}\rho(\pm n)n^{-s},\qquad
\mathfrak g_0(s)=\pi^{2s-1}\frac{\Ga(\frac{1-s}2)^2}{\Ga(\frac s2)^2},\qquad\mathfrak g_1(s)=\pi^{2s-1}\frac{\Ga(\frac{2-s}2)^2}{\Ga(\frac{1+s}2)^2}.
\]

\begin{theorem}\label{thm:FE}
$L_Q(s,h/k)$ continues meromorphically to $\C$. Its only possible singularity is a simple pole at $s=1$, with residue
$\log\eps_0\,\GQ(k,h)/(k^2\sqrt\Delta)$. Moreover
\begin{equation}\label{eq:FE}
L_Q(s,h/k)=\tfrac12\Delta^{\frac12-s}k^{-2s}\Bigl\{(\mathfrak g_0+\mathfrak g_1)(s)D_+(1-s)+(\mathfrak g_0-\mathfrak g_1)(s)D_-(1-s)\Bigr\}\qquad(\sigma<0).
\end{equation}
\end{theorem}

\begin{proof}
Let $m$ be the order of $\gamma_0$ modulo $k$. Put $\Ga_1=\langle\gamma_0^m\rangle$ and $\eps_1=\eps_0^m$; then $\Ga_1$ fixes every
class $a+k\Z^2$, and $[\Ga:\Ga_1]=2m$. For $p,q\in\{0,1\}$ define
\[
Z_{p,q}(s;a)=\sum_{\substack{v\in\Ga_1\backslash(a+k\Z^2)\\v\ne0}}\frac{\sgn\lambda_1(v)^p\sgn\lambda_2(v)^q}{|Q(v)|^s},
\]
\[
Z^*_{p,q}(s;a)=\sum_{\substack{w\in\Ga_1^*\backslash\Z^2\\w\ne0}}\e_k(a\cdot w)\frac{\sgn\mu_1(w)^p\sgn\mu_2(w)^q}{|\mu_1(w)\mu_2(w)|^s},
\]
and $G_{p,q}(s)=\pi^{-s-\frac{p+q}2}\Ga(\frac{s+p}2)\Ga(\frac{s+q}2)$.

\emph{Step 1: integral representation.} Let
\[
\theta(t_1,t_2)=\sum_{v\in a+k\Z^2}\lambda_1(v)^p\lambda_2(v)^q\e^{-\pi(t_1\lambda_1(v)^2+t_2\lambda_2(v)^2)},
\]
and let $\theta^\circ$ omit $v=0$. The Mellin transform of $x^p\e^{-\pi tx^2}$ gives, for $v\ne0$ and $\lambda_i=\lambda_i(v)$,
\[
\int_0^\infty\!\!\int_0^\infty\lambda_1^p\lambda_2^q\e^{-\pi(t_1\lambda_1^2+t_2\lambda_2^2)}t_1^{\frac{s+p}2}t_2^{\frac{s+q}2}\frac{dt_1dt_2}{t_1t_2}
=G_{p,q}(s)\frac{\sgn\lambda_1^p\sgn\lambda_2^q}{|\lambda_1\lambda_2|^{s}}.
\]
Substitute $t_1=tu$, $t_2=t/u$. Since $\gamma_0^m$ multiplies $u$ by $\eps_1^2$, unfolding gives
\begin{equation}\label{eq:hecke}
G_{p,q}(s)Z_{p,q}(s;a)=\int_0^\infty\mathcal F^\circ(t)\,t^{s+\frac{p+q}2}\frac{dt}t,\qquad\mathcal F(t)=2\int_1^{\eps_1^2}u^{\frac{p-q}2}\theta(tu,t/u)\frac{du}u ,
\end{equation}
where $\mathcal F^\circ$ is formed in the same way from $\theta^\circ$.

\emph{Step 2: Poisson summation.} We use $|\det\Theta|=\sqrt\Delta$ and
\[
\sum_{v\in a+k\Z^2}F(v)=k^{-2}\sum_w\e_k(a\cdot w)\widehat F(w/k),\qquad
\int_\R y^p\e^{-\pi ty^2}\e(-y\mu)\,dy=(-i)^pt^{-\frac12-p}\mu^p\e^{-\pi\mu^2/t}.
\]
With $t'=1/(k^2t)$, $u'=1/u$ and $\theta^*(T_1,T_2)=\sum_w\e_k(a\cdot w)\mu_1^p\mu_2^q\e^{-\pi(T_1\mu_1^2+T_2\mu_2^2)}$ this gives
\[
\theta(tu,t/u)=(-i)^{p+q}k^{-2-p-q}\Delta^{-1/2}t^{-1-p-q}u^{q-p}\theta^*(t'u',t'/u').
\]
Hence $\mathcal F(t)=(-i)^{p+q}k^{-2-p-q}\Delta^{-1/2}t^{-1-p-q}\mathcal F^*(t')$, where $\mathcal F^*$ is formed from $\theta^*$ and unfolded with
respect to $\Ga_1^*$.

\emph{Step 3: splitting at $t=1/k$.} As in Riemann's proof this gives the meromorphic continuation and
\begin{equation}\label{eq:FEpq}
G_{p,q}(s)Z_{p,q}(s;a)=(-i)^{p+q}\Delta^{-1/2}k^{-2s}G_{p,q}(1-s)Z^*_{p,q}(1-s;a).
\end{equation}
For $p=q=0$ the terms $w=0$ and $v=0$ contribute the constants $4\log\eps_1$ to $\mathcal F^*$ and $4\log\eps_1\mathbf 1_{a\equiv0}$ to
$\Phi$. They produce
\[
\frac{4\log\eps_1}{\sqrt\Delta}\cdot\frac{k^{-1-s}}{s-1}-4\log\eps_1\mathbf 1_{a\equiv0}\cdot\frac{k^{-s}}s .
\]
So the only pole of $Z_{0,0}$ is at $s=1$, with residue $4\log\eps_1/(k^2\sqrt\Delta)$ (as $G_{0,0}(1)=1$); the pole at $s=0$ is
cancelled by $G_{0,0}$. For $(p,q)\ne(0,0)$, $Z_{p,q}$ is entire.

\emph{Step 4: assembling.} Since $\mathbf 1_{Q>0}=\frac12(1+\sgn\lambda_1\sgn\lambda_2)$ and $\Ga$ acts freely,
\[
L_Q(s,h/k)=\frac1{4m}\sum_{a\bmod k}\e_k(hQ(a))\bigl(Z_{0,0}+Z_{1,1}\bigr)(s;a).
\]
Apply \eqref{eq:FEpq}, noting $G_{p,p}(1-s)/G_{p,p}(s)=\mathfrak g_p(s)$ and $(-i)^2=-1$. By \eqref{eq:mumu},
$|\mu_1\mu_2|=|Q^*|/\Delta$ and $\sgn\mu_1\sgn\mu_2=-\sgn Q^*$. Further $\sum_{a\bmod k}\e_k(hQ(a))\e_k(a\cdot w)=\GQ(k,h,w)$, and by
Lemma~\ref{lem:gauss}(a) and freeness $\sum_{\Ga_1^*\backslash}=2m\sum_{\Ga^*\backslash}$ for $\Ga^*$-invariant summands. Collecting
terms,
\[
L_Q(s,h/k)=\tfrac12\Delta^{\frac12-s}k^{-2s}\Bigl\{\mathfrak g_0(s)\sum_{n\ne0}\rho(n)|n|^{s-1}+\mathfrak g_1(s)\sum_{n\ne0}\sgn(n)\rho(n)|n|^{s-1}\Bigr\},
\]
which is \eqref{eq:FE}. The residue is $\frac1{4m}\GQ(k,h)\frac{4\log\eps_1}{k^2\sqrt\Delta}=\frac{\log\eps_0\GQ(k,h)}{k^2\sqrt\Delta}$.
\end{proof}

\begin{corollary}\label{cor:selfdual}
For $k=1$, $\rho_Q(n;0)=r_Q(n)$. The function $\Lambda$ of \eqref{eq:FEintro} and
\[
\Lambda^-(s)=\Bigl(\frac{\sqrt\Delta}\pi\Bigr)^s\Ga\Bigl(\frac{s+1}2\Bigr)^2\Xi^-_Q(s)
\]
satisfy $\Lambda(s)=\Lambda(1-s)$ and $\Lambda^-(s)=\Lambda^-(1-s)$. $\zQ$ has a simple pole at $s=1$ with residue
$2\log\eps_0/\sqrt\Delta$, and $\Xi^-_Q$ is entire.
\end{corollary}
\begin{proof}
Add, respectively subtract, \eqref{eq:FE} for $Q$ and for $-Q$, using $r_{-Q}(n)=r_Q(-n)$.
\end{proof}

\begin{remark}[Narrow classes]\label{rem:narrow}
$Z_Q=\frac12(\zQ+\Xi^-_Q)$, and $\mathfrak g_0(s)/\mathfrak g_1(s)=\tan^2(\pi s/2)\to-1$ exponentially fast on $\sigma=\frac12$. Since
$\Lambda$ and $\Lambda^-$ are real on the critical line, $Z_Q(\frac12+it)$ equals a common phase times $W_+(t)-iW_-(t)$, up to
$O(\e^{-\pi t})$, with $W_\pm$ real. A real zero requires $W_+$ and $W_-$ to vanish together. So Hardy's method gives
nothing for $Z_Q$ unless $\Xi^-_Q\equiv0$, e.g.\ when $Q\sim-Q$.
\end{remark}

\begin{lemma}\label{lem:mellin}
For $\Re w>-1$,
\[
(\mathfrak g_0+\mathfrak g_1)(-w)=\pi^{-2}(2\pi)^{-2w}\Ga(1+w)^2=\int_0^\infty8K_0(4\pi\sqrt y)\,y^{w}\,dy,
\]
and for $-1<\Re w<-\frac14$,
\[
(\mathfrak g_0-\mathfrak g_1)(-w)=-\pi^{-2}(2\pi)^{-2w}\Ga(1+w)^2\cos\pi w=\int_0^\infty\bigl(-4\pi Y_0(4\pi\sqrt y)\bigr)\,y^{w}\,dy,
\]
the last integral being convergent but not absolutely convergent.
\end{lemma}
\begin{proof}
By the reflection and duplication formulas,
\[
\mathfrak g_0(-w)=\pi^{-2}(2\pi)^{-2w}\Ga(1+w)^2\sin^2\frac{\pi w}2,\qquad\mathfrak g_1(-w)=\pi^{-2}(2\pi)^{-2w}\Ga(1+w)^2\cos^2\frac{\pi w}2.
\]
Now substitute $x=4\pi\sqrt y$ and use
\[
\int_0^\infty K_0(x)x^{z-1}dx=2^{z-2}\Ga(\tfrac z2)^2,\qquad
\int_0^\infty Y_0(x)x^{z-1}dx=-\frac{2^{z-1}}\pi\Ga(\tfrac z2)^2\cos\frac{\pi z}2\quad(0<\Re z<\tfrac32).\qedhere
\]
\end{proof}

\begin{proof}[Proof of Theorem~\ref{thm:voronoi}]
As in the proof of Proposition~\ref{prop:Vdef}, shift $\frac1{2\pi i}\int_{(2)}L_Q(s,h/k)\tilde f(s)\,ds$ to $\sigma=-\frac18$, picking
up the residue at $s=1$; on $\sigma=-\frac18$ we have $L_Q\ll(1+|s|)^{5/4}$ by \eqref{eq:FE} and Stirling, and $L_Q$ has finite
order by \eqref{eq:hecke}. Insert \eqref{eq:FE} and integrate term by term. The term $\rho(\pm n)$ contributes
\[
\frac{\sqrt\Delta}{2n}\cdot\frac1{2\pi i}\int(\mathfrak g_0\pm\mathfrak g_1)(s)\,y^s\tilde f(s)\,ds,\qquad y=\frac n{\Delta k^2}.
\]
For the minus sign, first move to $\sigma\in(\frac14,1)$; this is allowed since the poles of $\mathfrak g_0$ are at $1,3,5,\dots$ and
those of $\mathfrak g_1$ at $2,4,\dots$. Then apply Lemma~\ref{lem:mellin} and Mellin--Parseval. Absolute convergence follows from
$\int fY_0(c\sqrt x)\,dx\ll_fc^{-7/2}$ (three integrations by parts, as before) and $\rho(n)\ll_kn^\eps$.
\end{proof}

\begin{corollary}\label{cor:epstein}
Let $\omega\in\{\pm1\}$ and $\rho^\omega(m)=\rho_Q(m;h/k)+\omega\rho_Q(-m;-h/k)$. Then $\sum_{n\ge1}(r_Q(n)+\omega r_Q(-n))\e(hn/k)f(n)$
equals the right-hand side of Theorem~\ref{thm:voronoi} with $\GQ(k,h)$ replaced by $\GQ(k,h)+\omega\GQ(k,-h)$ and
$\rho(\pm n)$ by $\rho^\omega(\pm n)$.
\end{corollary}
\begin{proof}
Apply Theorem~\ref{thm:voronoi} to $-Q$, noting $r_{-Q}(n)=r_Q(-n)$, $(-Q)^*=-Q^*$, $G_{-Q}(k,h,w)=\GQ(k,-h,w)$ and
$\Ga_{-Q}=\Ga$, so that $\rho_{-Q}(m;h/k)=\rho_Q(-m;-h/k)$.
\end{proof}

\subsection{The twisted formulas in unified form}
For the application we need the formulas with the twist $\e(-hn/k)$, $(k,2\Delta)=1$, and with the Bessel functions split
into Hankel functions. Let $H_0^{(+1)}=H_0^{(1)}$ and $H_0^{(-1)}=H_0^{(2)}$, so that
\[
J_0=\tfrac12\bigl(H_0^{(+1)}+H_0^{(-1)}\bigr),\qquad-2\pi Y_0=\pi\bigl(iH_0^{(+1)}-iH_0^{(-1)}\bigr).
\]
Recall the data $\Dc$, $P$, $\beta$, $\theta_\vs$ of Table~\ref{tab:cases}, and put
\begin{equation}\label{eq:qdef}
q=k\sqrt\Delta .
\end{equation}

\begin{corollary}\label{cor:unified}
Let $(k,2\Delta)=1$, $(h,k)=1$, and let $f\in C^3$ have compact support in $(0,\infty)$. Then
\begin{multline}\label{eq:unified}
\sum_{n\ge1}c(n)\e\Bigl(-\frac{hn}k\Bigr)f(n)=\mathfrak m\int_0^\infty f\\
+\frac\pi q\sum_{w\in\Dc}\beta(w)\,\e_k\bigl(\hb\Db P(w)\bigr)\sum_{\vs=\pm1}\theta_\vs\int_0^\infty f(x)H_0^{(\vs)}\Bigl(\frac{4\pi\sqrt{P(w)x}}q\Bigr)dx
+\mathfrak k(f),
\end{multline}
where $|\beta(w)|=1$, and:
\begin{itemize}
\item definite case: $\mathfrak m=2\pi\GQ(k,-h)/(k^2\sqrt\Delta)$ and $\mathfrak k(f)=0$;
\item indefinite case: $\mathfrak m=\log\eps_0\,(\GQ(k,-h)+\omega\GQ(k,h))/(k^2\sqrt\Delta)$ and
\[
\mathfrak k(f)=\frac4{k^2\sqrt\Delta}\sum_{n\ge1}\rho^\omega(n)\int_0^\infty f(x)K_0\Bigl(\frac{4\pi\sqrt{nx}}q\Bigr)dx,
\qquad|\rho^\omega(n)|\le k\,\#\{w\in\Dc:P(w)=n\},
\]
with $\rho^\omega$ formed as in Corollary~\ref{cor:epstein} with $-h$ in place of $h$.
\end{itemize}
In both cases $|\mathfrak m|\ll1/k$.
\end{corollary}
\begin{proof}
\emph{Definite case.} Apply Proposition~\ref{prop:Vdef} with $-h$ in place of $h$. By Lemma~\ref{lem:gauss}(c) with
$\det A=\Delta$, $\GQ(k,-h,w)=\GQ(k,-h)\e_k(\hb\Db Q^*(w))=k\beta\,\e_k(\hb\Db P(w))$. Now write $J_0$ as above.

\emph{Indefinite case.} Apply Corollary~\ref{cor:epstein} with $-h$ in place of $h$. Then
\[
\rho^\omega(-n)=\rho_Q(-n;-h/k)+\omega\rho_Q(n;h/k)=\sum_{Q^*(w)=-n}\GQ(k,-h,w)+\omega\sum_{Q^*(w)=n}\GQ(k,h,w),
\]
the sums running over $\Ga^*$-orbits. By Lemma~\ref{lem:gauss}(c) with $\det A=-\Delta$,
\[
\GQ(k,-h,w)=\GQ(k,-h)\e_k(-\hb\Db Q^*(w)),\qquad\GQ(k,h,w)=\GQ(k,h)\e_k(\hb\Db Q^*(w)),
\]
so both sums carry the factor $\e_k(\hb\Db n)$, and $\rho^\omega(-n)=k\,\e_k(\hb\Db n)\sum_{w\in\Dc,P(w)=n}\beta(w)$. Inserting this
and $-2\pi Y_0=\pi\sum_\vs i\vs H_0^{(\vs)}$ gives the $Y_0$-part of \eqref{eq:unified}. The bound for $\rho^\omega(n)$ follows from
$|\GQ(k,\pm h,w)|=k$. Finally $|\GQ(k,\pm h)|=k$ gives $|\beta(w)|=1$ and $|\mathfrak m|\ll1/k$.
\end{proof}

\section{Oscillatory integrals}\label{sec:osc}

The Voronoi formulas of \S\ref{sec:FEV} express $\Pc(t)$ through Hankel transforms of the weight $\eta(x)x^{-1/2-it}\e(hx/k)$.
In this section we evaluate these transforms asymptotically (Lemma~\ref{lem:statphase}). We use only the classical
asymptotic expansion of the Hankel functions and two lemmas of Blomer, Khan and Young on oscillatory integrals.

\begin{fact}[Hankel functions; {\cite[\S\S6.12, 7.2, 7.23]{Wat}}]\label{fact:hankel}
For real $z>0$ and $\vs=\pm1$,
\[
H_0^{(\vs)}(z)=\sqrt{\frac2{\pi z}}\,\e^{i\vs(z-\pi/4)}\bigl(1+w_\vs(z)\bigr),\qquad
w_{+1}(z)=\frac1{\sqrt\pi}\int_0^\infty\e^{-u}u^{-1/2}\Bigl(\bigl(1+\tfrac{iu}{2z}\bigr)^{-1/2}-1\Bigr)du,
\]
and $w_{-1}=\overline{w_{+1}}$. Consequently $z^jw_\vs^{(j)}(z)\ll_jz^{-1}$ for $z\ge1$ and $j\ge0$. Moreover
$0<K_0(z)\le\sqrt{\pi/(2z)}\,\e^{-z}$.
\end{fact}

The bound for $w_\vs$ follows by differentiating under the integral: with $v=u/(2z)$ and $F(v)=(1+iv)^{-1/2}-1$, $z^j\partial_z^j$
is a combination of $v^iF^{(i)}(v)$, $i\le j$, each of which is $\ll\min(v,1)$ because $|1+iv|\ge1$.

\begin{fact}[First derivative test; {\cite[Prop.~8.1]{BKY}}]\label{fact:BKY1}
Let $Y\ge1$, $X,Q,U,R>0$, and let $w$ be smooth with support in $[\alpha,\beta]$, with $w^{(j)}\ll_jXU^{-j}$. Let $h$ be smooth on
$[\alpha,\beta]$ with $|h'|\ge R$ and $h^{(j)}\ll_jYQ^{-j}$ for $j\ge2$. Then
\[
\int w(x)\e^{ih(x)}dx\ll_A(\beta-\alpha)X\Bigl[\Bigl(\frac{QR}{\sqrt Y}\Bigr)^{-A}+(RU)^{-A}\Bigr].
\]
\end{fact}

\begin{fact}[Stationary phase; {\cite[Prop.~8.2]{BKY}}]\label{fact:BKY2}
Let $0<\delta<\frac1{10}$, $X,Y,V,V_1,Q>0$, $Z=Q+X+Y+V_1+1$, and assume $Y\ge Z^{3\delta}$ and $V_1\ge V\ge QZ^{\delta/2}Y^{-1/2}$.
Let $w$ be smooth with support in an interval $J'$ of length $V_1$, with $w^{(j)}\ll_jXV^{-j}$. Let $h$ be smooth on $J'$ with a
unique $x_0\in J'$ such that $h'(x_0)=0$, and with $h''\gg YQ^{-2}$ and $h^{(j)}\ll_jYQ^{-j}$ on $J'$ for $j\ge1$. Then
\[
\int w(x)\e^{ih(x)}dx=\frac{\e^{ih(x_0)}}{\sqrt{h''(x_0)}}\sum_{j\le3A/\delta}p_j(x_0)+O_{A,\delta}(Z^{-A}),\qquad
p_0(x_0)=\sqrt{2\pi}\,\e^{i\pi/4}w(x_0),
\]
with $p_j(x_0)\ll X\bigl((V^2Y/Q^2)^{-j}+Y^{-j/3}\bigr)$; moreover $p_j(x_0)$ depends only on the derivatives of $w$ and $h$ at
$x_0$, and vanishes if $w$ vanishes near $x_0$.
\end{fact}

\subsection{The Bessel transform of the critical weight}
Fix $\eta_0\in C^\infty(\R)$ with $0\le\eta_0\le1$, $\eta_0=1$ on $[-\frac12,\frac12]$, $\supp\eta_0\subset[-1,1]$, and $\eta_0$
non-decreasing on $(-\infty,0]$ and non-increasing on $[0,\infty)$.

\begin{lemma}\label{lem:statphase}
Let $T$ be large, $T^{1/2+\delta}\le K\le T^{1-\delta}$, and let $x_*,t,r,q$ satisfy
\begin{gather*}
c_*T\le x_*\le C_*T,\qquad\tfrac12T\le t\le2T,\qquad c_r\le r\le C_r,\qquad1\le q\le T^{1/2-\delta},\\
|x_0-x_*|\le\frac K2,\qquad\text{where }x_0=\frac\tau r,\ \ \tau=\frac t{2\pi}.
\end{gather*}
Let $\eta(x)=\eta_0((x-x_*)/K)$, $J=[x_*-K,x_*+K]$, $g(x)=rx-\tau\log x$, and
\[
N_0=\frac{r^2q^2K^2}{x_*} .
\]
For $u>0$ and $\vs=\pm1$ define
\[
s_u=\Bigl(\frac{\pi u}{2q^2rt}\Bigr)^{1/2},\qquad x_\vs(u)=x_0\,\e^{-2\vs\arsinh s_u},\qquad z_\vs(u)=\frac{4\pi\sqrt{u\,x_\vs(u)}}q,
\qquad\Phi(u)=\frac t\pi\,\phi(s_u^2),
\]
\[
\mathcal M_\vs(u)=\Bigl(\frac q{2\pi^2}\Bigr)^{1/2}\e^{i(1-\vs)\pi/4}\,u^{-1/4}\,\eta(x_\vs(u))\,x_\vs(u)^{-3/4}\bigl(1+w_\vs(z_\vs(u))\bigr)\,g_\vs''(x_\vs(u))^{-1/2},
\]
where $g_\vs''(x_\vs(u))=\frac{r^2}\tau\,\e^{4\vs\arsinh s_u}\bigl(1-\vs s_u\e^{-\vs\arsinh s_u}\bigr)$. For $n>0$ let
\[
\mathcal J_\vs(n)=\int_0^\infty\eta(x)x^{-1/2}\e(g(x))\,H_0^{(\vs)}\Bigl(\frac{4\pi\sqrt{nx}}q\Bigr)dx .
\]
Then the following hold, with implied constants depending on $\delta,c_*,C_*,c_r,C_r,\eta_0$ and $A$.
\begin{enumerate}
\item[(i)] $\displaystyle\mathcal J_\vs(n)=\e(\tau-\tau\log x_0)\,\e\Bigl(-\frac n{2q^2r}+\vs\Phi(n)\Bigr)\mathcal M_\vs(n)+\mathcal E_\vs(n)$, where
\[
\mathcal E_\vs(n)\ll q^{1/2}n^{-1/4}T^{-1/4}\Bigl(\frac T{K^2}+T^{-1/3}\Bigr)\quad(n\le20N_0),\qquad
\mathcal E_\vs(n)\ll_AT^{-A}\Bigl(\frac n{N_0}\Bigr)^{-A}\quad(n>20N_0).
\]
\item[(ii)] $\mathcal M_\vs(u)=0$ unless $x_\vs(u)\in J$, and this implies $u\le5N_0$. Moreover, for all $u,N\ge1$,
\[
|\mathcal M_\vs(u)|\ll q^{1/2}u^{-1/4}T^{-1/4},\qquad\Var_{[N,2N]}\mathcal M_\vs\ll q^{1/2}N^{-1/4}T^{-1/4}.
\]
\item[(iii)] $\displaystyle\int_0^\infty\eta(x)x^{-1/2}\e(g(x))K_0\Bigl(\frac{4\pi\sqrt{nx}}q\Bigr)dx\ll KT^{-1/2}\exp\Bigl(-c\frac{\sqrt{nT}}q\Bigr)$ for some $c>0$.
\item[(iv)] $\displaystyle\int_0^\infty\eta(x)x^{-1/2}\e(g(x))\,dx=\e^{i\pi/4}r^{-1/2}\e(\tau-\tau\log x_0)+O\Bigl(\frac T{K^2}+T^{-1/3}\Bigr)$; in particular
it is $\ll1$.
\end{enumerate}
\end{lemma}

\begin{figure}[t]
\centering
\input{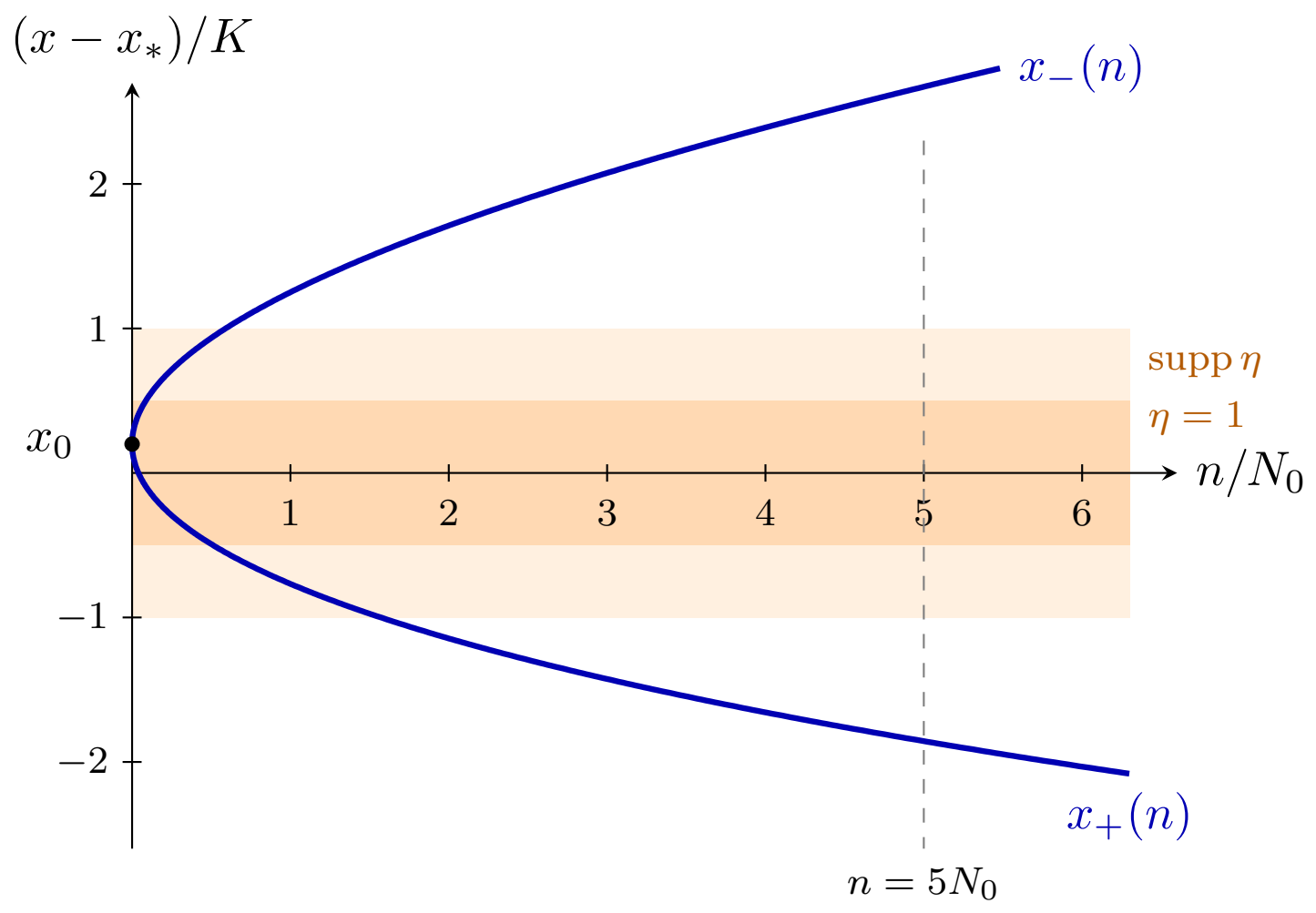}
\caption{The stationary points $x_\pm(n)$ of Lemma~\ref{lem:statphase}, given exactly by $n=r^2q^2(x_0-x)^2/x$ (here
$x_*=12K$ and $x_0-x_*=0.2K$). The term $n$ of the dual sum contributes only if $x_\vs(n)$ lies in $\supp\eta$ (light), which
forces $n\le5N_0$; for larger $n$ the integral is negligible by Fact~\ref{fact:BKY1}.}
\label{fig:statpoints}
\end{figure}

\begin{proof}
\emph{Reduction to an exponential integral.} By Fact~\ref{fact:hankel},
\begin{equation}\label{eq:Jred}
\mathcal J_\vs(n)=\Bigl(\frac q{2\pi^2}\Bigr)^{1/2}\e^{-i\vs\pi/4}n^{-1/4}\int W_\vs(x)\e^{ih_\vs(x)}dx,
\end{equation}
with
\[
W_\vs(x)=\eta(x)x^{-3/4}\Bigl(1+w_\vs\Bigl(\frac{4\pi\sqrt{nx}}q\Bigr)\Bigr),\qquad h_\vs=2\pi g_\vs,\qquad g_\vs(x)=rx-\tau\log x+\frac{2\vs\sqrt{nx}}q .
\]
On $J$ we have $4\pi\sqrt{nx}/q\gg T^{\delta}$, since $q\le T^{1/2-\delta}$; so $|w_\vs|\le\frac12$ there, and
$W_\vs^{(j)}\ll_jT^{-3/4}K^{-j}$.

\emph{The stationary point.} Since
\[
x\,g_\vs'(x)=ry^2+\vs\frac{\sqrt n}qy-\tau\qquad(y=\sqrt x),
\]
$g_\vs$ has exactly one stationary point on $(0,\infty)$, at
\[
\sqrt{x_\vs}=\sqrt{\tau/r}\,\bigl(\sqrt{1+s_n^2}-\vs s_n\bigr)=\sqrt{x_0}\,\e^{-\vs\arsinh s_n},
\]
since $\sqrt n/(2q\sqrt{r\tau})=s_n$. From $g_\vs'(x_\vs)=0$ we get $\vs\sqrt n/(q\sqrt{x_\vs})=\tau/x_\vs-r=r(x_0-x_\vs)/x_\vs$, hence
\begin{equation}\label{eq:nx}
n=\frac{r^2q^2(x_0-x_\vs)^2}{x_\vs}.
\end{equation}
Using $rx_\vs=\tau-\vs\sqrt n\sqrt{x_\vs}/q$ and $\sqrt n\sqrt{x_\vs}/q=2\tau s_n\e^{-\vs\arsinh s_n}$, we obtain
\[
g_\vs(x_\vs)=\tau-\tau\log x_0+2\vs\tau\bigl(\arsinh s_n+s_n\e^{-\vs\arsinh s_n}\bigr)
=\tau-\tau\log x_0-2\tau s_n^2+2\vs\tau\,\phi(s_n^2),
\]
because $s\e^{\mp\arsinh s}=s\sqrt{1+s^2}\mp s^2$. Here $2\tau s_n^2=n/(2q^2r)$ and $2\tau\phi(s_n^2)=\Phi(n)$. Similarly,
$g_\vs''(x)=x^{-2}\bigl(\tau-\frac\vs2\sqrt n\sqrt x/q\bigr)$ gives the stated value of $g''_\vs(x_\vs)$.

\emph{Case 1: $x_\vs\in J'=[x_*-2K,x_*+2K]$.} Then $|x_0-x_\vs|\le\frac52K$ and $x_\vs\ge x_*/2$, so $n\le12.5N_0$ by
\eqref{eq:nx}. Hence $\sqrt n/q\le4rK/\sqrt{x_*}$, and on $J'$ the term $2\vs\sqrt{nx}/q$ and its derivatives satisfy
$\partial_x^j(\sqrt{nx}/q)\ll KT^{-j}$. Since $K\le T^{1-\delta}$, this gives on $J'$
\[
h_\vs''=2\pi\Bigl(\frac\tau{x^2}+O\Bigl(\frac K{T^2}\Bigr)\Bigr)\asymp\frac1T,\qquad h_\vs^{(j)}\ll_jT^{1-j}.
\]
In particular $h_\vs'$ is increasing on $J'$, and $x_\vs$ is its unique zero there. Apply Fact~\ref{fact:BKY2} with
\[
X\asymp T^{-3/4},\quad V=K,\quad V_1=4K,\quad Y=Q=T .
\]
Its hypotheses hold since $Z\ll T$ and $K\ge T^{1/2+\delta}$. The leading term is
\[
\frac{\e^{ih_\vs(x_\vs)}}{\sqrt{h_\vs''(x_\vs)}}\sqrt{2\pi}\,\e^{i\pi/4}W_\vs(x_\vs)=\e(g_\vs(x_\vs))\,\e^{i\pi/4}\,\frac{W_\vs(x_\vs)}{g_\vs''(x_\vs)^{1/2}},
\]
and the other terms are $\ll X/\sqrt{h''}\cdot(T/K^2+T^{-1/3})\ll T^{-1/4}(T/K^2+T^{-1/3})$. Inserting this into \eqref{eq:Jred} gives (i)
for these $n$.

\emph{Case 2: $x_\vs\notin J'$.} Then $\eta(x_\vs)=0$, so $\mathcal M_\vs(n)=0$ and $\mathcal E_\vs(n)=\mathcal J_\vs(n)$. Let $x\in J$ and
$y=\sqrt x$, $y_\vs=\sqrt{x_\vs}$. Then $|x-x_\vs|\ge K$, and
\[
x\,g'_\vs(x)=r(y-y_\vs)\Bigl(y+\frac\tau{ry_\vs}\Bigr),\qquad\text{so}\qquad|h'_\vs(x)|\ge\frac{2\pi r|y-y_\vs|}{\sqrt x}\ge\frac{2\pi rK}{3x}.
\]
(If $x_\vs\le4x$, then $|y-y_\vs|=|x-x_\vs|/(y+y_\vs)\ge K/(3\sqrt x)$; if $x_\vs>4x$, then $|y-y_\vs|\ge\sqrt x\ge K/(3\sqrt x)$.) If
$n>32N_0$, then also $\sqrt n\,y/q\ge4rK\ge2|rx-\tau|$ on $J$, hence $|h_\vs'(x)|\ge\pi\sqrt n/(q\sqrt x)$. On $J$,
$h_\vs^{(j)}\ll(T+\sqrt{nT}/q)T^{-j}$ for $j\ge2$. Apply Fact~\ref{fact:BKY1} with $U=K$, $X\asymp T^{-3/4}$, $Q=T$,
$Y=T+\sqrt{nT}/q$:
\begin{itemize}
\item for $n\le32N_0$ take $R\asymp K/T$; then $QR/\sqrt Y\gg K/\sqrt T\ge T^\delta$ and $RU\gg K^2/T\ge T^{2\delta}$;
\item for $n>32N_0$ take $R\asymp\sqrt n/(q\sqrt T)$. Since $\sqrt{N_0}/q\asymp K/\sqrt T\ge T^\delta$, we get
$QR/\sqrt Y\gg\min\bigl(\sqrt n/q,(\sqrt{nT}/q)^{1/2}\bigr)\gg T^{\delta}(n/N_0)^{1/4}$ and $RU\gg T^{2\delta}(n/N_0)^{1/2}$.
\end{itemize}
In both cases $\mathcal J_\vs(n)\ll_AT^{-A}(1+n/N_0)^{-A}$. Since Case~1 forces $n\le12.5N_0$, every $n>20N_0$ is in Case~2. This
proves (i).

(ii) If $x_\vs(u)\in J$, then $|x_0-x_\vs|\le\frac32K$ and $x_\vs\ge x_*/2$, so $u\le4.5N_0$ by \eqref{eq:nx}. For $u\le5N_0$ we
have $s_u^2\ll K^2/T^2\le T^{-2\delta}$, so $x_\vs(u)\asymp T$, $g''_\vs(x_\vs(u))\asymp1/T$ and $|1+w_\vs|\asymp1$; this gives the bound for
$|\mathcal M_\vs|$ (note $z_\vs(u)\gg T^\delta$ for $u\ge1$). For the variation write $\mathcal M_\vs(u)=c\,\eta(x_\vs(u))F(u)$. Each factor $F_i$ of $F$ (namely $u^{-1/4}$,
$x_\vs^{-3/4}$, $1+w_\vs(z_\vs)$ and $g''_\vs(x_\vs)^{-1/2}$) satisfies $|uF_i'(u)|\ll|F_i(u)|$, because $u\frac d{du}\log x_\vs=-\vs s_u/\sqrt{1+s_u^2}$,
$u\frac d{du}s_u=\frac12s_u$, and $zw_\vs'(z)\ll z^{-1}$. Hence $\Var_{[N,2N]}F\ll\sup_{[N,2N]}|F|$ for $1\le N\le5N_0$, and $\sup_{[N,2N]}|F|\ll N^{-1/4}T^{-1/4}$ there. Since
$u\mapsto x_\vs(u)$ is monotone and $\eta$ is unimodal, $\Var(\eta\circ x_\vs)\le2$. Therefore
\[
\Var_{[N,2N]}(\eta\circ x_\vs\cdot F)\le\sup_{[N,2N]}|F|\cdot\Var(\eta\circ x_\vs)+\Var_{[N,2N]}F\ll N^{-1/4}T^{-1/4},
\]
which gives the claim for $N\le5N_0$; for $N>5N_0$, $\mathcal M_\vs$ vanishes on $[N,2N]$.

(iii) On $J$, $4\pi\sqrt{nx}/q\ge c\sqrt{nT}/q\ge2$, so Fact~\ref{fact:hankel} bounds the integrand by $T^{-1/2}\e^{-c\sqrt{nT}/q}$.

(iv) Apply Fact~\ref{fact:BKY2} to $w=\eta x^{-1/2}$ ($X\asymp T^{-1/2}$) and $h=2\pi g$, which has the stationary point $x_0$ with
$g''(x_0)=r^2/\tau$ and $g(x_0)=\tau-\tau\log x_0$. Since $\eta(x_0)=1$, the leading term is
$\e(g(x_0))\e^{i\pi/4}x_0^{-1/2}(\tau/r^2)^{1/2}=\e^{i\pi/4}r^{-1/2}\e(\tau-\tau\log x_0)$, and the other terms are
$\ll X(T/r^2)^{1/2}(T/K^2+T^{-1/3})$.
\end{proof}

\section{Transformation of the critical Dirichlet polynomial}\label{sec:transf}

In this section $Q$ is as in Theorem~\ref{thm:A}(a) or Theorem~\ref{thm:B}, and $Z$, $c(n)$, $\omega$ are as in \S\ref{subsec:language};
if $\omega=-1$ we assume $\Xi^-_Q\not\equiv0$.

\subsection{Zero detection}\label{subsec:zerodet}
Put
\begin{equation}\label{eq:xstar}
x_*=\frac{T\sqrt\Delta}{2\pi}.
\end{equation}

\begin{proposition}\label{prop:zerodet}
Let $T$ be large, $T^{2\eps}\le H\le T^{1/2}$, $K=T^{1+2\eps}/H$, and let $\eta:\R\to[0,1]$ satisfy $\eta=1$ on
$[x_*-\frac K2,x_*+\frac K2]$ and $\eta=0$ outside $[x_*-K,x_*+K]$. If
\[
\Pc(t):=\sum_n\eta(n)c(n)n^{-1/2-it}\ll(\log T)^{-2}\qquad(|t-T|\le H),
\]
then $Z$ has a zero $\frac12+i\gamma$ with $|\gamma-T|\le H$.
\end{proposition}

In the definite case this is \cite[\S2]{JS}, with the details in \cite[\S2]{BSS}; the argument uses only the stated
properties of $\eta$ (they are listed before \cite[(2.4)]{JS}). The approximate formula used there \cite[Lemma~2]{BSS}
(stated without proof in \cite[Lemma~1]{JS}) can be replaced by the following lemma, whose slightly weaker error term
is sufficient. We state it for both cases, since we need it for indefinite $Q$.

\begin{lemma}\label{lem:approx}
Let $R_Z$ be the residue of $Z$ at $s=1$, i.e.\ $R_Z=2\pi/\sqrt\Delta$ for $Z=\zeta_Q$, $R_Z=2\log\eps_0/\sqrt\Delta$ for $Z=\zQ$, and
$R_Z=0$ for $Z=\Xi^-_Q$. Let $0<\xi\le\frac14$ and $A>0$. For $t\ge2$, $s=\frac12+it$ and $t^2\le X\le t^A$,
\[
Z(s)=\sum_{n\le X}\frac{c(n)}{n^{s}}+\frac1{\log2}\sum_{X<n\le2X}\frac{c(n)\log(2X/n)}{n^{s}}
-\frac{R_Z}{\log2}\cdot\frac{(2X)^{1-s}-X^{1-s}}{(1-s)^2}+O_{\xi,A}\bigl(tX^{-1/2+\xi}\bigr).
\]
\end{lemma}
\begin{proof}
By the functional equation (Proposition~\ref{prop:FEdef} with $k=1$, respectively Corollary~\ref{cor:selfdual}) and
Stirling's formula, $Z(-\frac\xi2+i\tau)\ll(1+|\tau|)^{1+\xi}$, and trivially $Z(1+\frac\xi2+i\tau)\ll1$. Since $(s-1)Z(s)$
is entire of finite order, the Phragm\'en--Lindel\"of principle gives
\begin{equation}\label{eq:convex}
Z(\sigma+i\tau)\ll(1+|\tau|)^{1+\frac\xi2-\sigma}\qquad(-\tfrac\xi2\le\sigma\le1+\tfrac\xi2,\ |\sigma+i\tau-1|\ge\tfrac12).
\end{equation}
Consider $\frac1{2\pi i}\int_{(1)}Z(s+w)X^w(2^w-1)w^{-2}\,dw$. Expanding $Z(s+w)$ and using
$\frac1{2\pi i}\int_{(1)}y^ww^{-2}dw=\log^+y$, it equals $\log2$ times the two sums. Shift the contour to $\Re w=-\frac12+\xi$,
i.e.\ $\Re(s+w)=\xi$; this is justified by \eqref{eq:convex}, since the exponent $1+\frac\xi2-\Re(s+w)$ is less than $1$ on and to
the right of the new line, while $|w|^{-2}$ decays. The residue at $w=0$ is $\log2\cdot Z(s)$, and the residue at $w=1-s$ is
$R_ZX^{1-s}(2^{1-s}-1)(1-s)^{-2}$. By \eqref{eq:convex}, the integral over the new line is
\[
\ll X^{-\frac12+\xi}\int_{-\infty}^\infty\frac{(1+|t+v|)^{1-\xi/2}}{(1+|v|)^2}\,dv\ll_\xi tX^{-\frac12+\xi}.
\]
Dividing by $\log2$ gives the lemma.
\end{proof}

\begin{proof}[Proof of Proposition~\ref{prop:zerodet} for indefinite $Q$]
Let $W(t)=\e^{\pi t/2}\Lambda(\frac12+it)$ if $\omega=1$ and $W(t)=\e^{\pi t/2}\Lambda^-(\frac12+it)$ if $\omega=-1$. By
Corollary~\ref{cor:selfdual}, $W$ is real, and by Stirling $|W(t)|\asymp t^{-\omega/2}|Z(\frac12+it)|$. Put $\omega'=(1-\omega)/2\in\{0,1\}$. The phase of
$(\sqrt\Delta/\pi)^{it}\Ga(\frac14+\frac{\omega'}2+\frac{it}2)^2$ has derivative $\log\frac{t\sqrt\Delta}{2\pi}+O(t^{-1})$, so
the resonance point is again \eqref{eq:xstar}. The argument of \cite[\S2]{BSS} now applies with three changes:
\begin{enumerate}
\item the lower bound for $I=\int_{-H}^HW(T+u)\e^{-(u/H_0)^2}du$, $H_0=HT^{-\eps}$, is $|I|\gg H_0T^{-\omega/2}$; it comes from
Ramachandra's theorem \cite[Lemma~1]{BSS} applied to $Z$, which is not identically zero;
\item the approximate formula \cite[Lemma~2]{BSS} is replaced by Lemma~\ref{lem:approx} with $X=T^3$ and $\xi=\frac1{12}$, whose
error term contributes $\ll H_0T^{-\omega/2}T^{-1/4}$ to $I$;
\item in the estimates of the non-critical terms, $\Ga(\frac12+i(T+u))$ is replaced by
$(\sqrt\Delta/\pi)^{i(T+u)}\Ga(\frac14+\frac{\omega'}2+\frac i2(T+u))^2$, whose $u$-phase has derivative
$\log\frac{T\sqrt\Delta}{2\pi}+O(T^{-1}(1+|u|))$; so the contour shifts $u\mapsto u\mp iH_0$ of \cite{BSS} give the same bounds
$\ll H\e^{-T^\eps}$.
\end{enumerate}
All upper bounds carry the same factor $T^{-\omega/2}$ as the lower bound, and the argument concludes as in \cite{BSS}.
\end{proof}

\subsection{Choice of the twist}\label{subsec:twist}
From now on
\begin{equation}\label{eq:Krange}
T^{1/2+\delta}\le K\le T^{3/5},\qquad HK=T^{1+2\eps},\qquad0<\eps\le\delta/4 .
\end{equation}

\begin{lemma}\label{lem:hk}
Let $\Delta$ be a positive non-square. There are $c_1,c_2,c_3,c_4>0$, depending on $\Delta$, such that for every $Y\ge c_1$ there
exist coprime $h,k$ with
\[
(k,2\Delta)=1,\qquad c_2Y\le k\le c_3Y,\qquad\Bigl|\frac1{\sqrt\Delta}-\frac hk\Bigr|\le\frac{c_4}{k^2}.
\]
\end{lemma}
\begin{proof}
The partial quotients of $\theta=1/\sqrt\Delta$ are bounded, say by $A$, so $|\theta-p/q|\ge((A+2)q^2)^{-1}$ for all $p/q$. Let
$p_j/q_j$ be the convergents with $q_j\le Y<q_{j+1}\le(A+2)q_j$, and let $P_0=\prod_{p\mid2\Delta}p$. Since $(q_j,q_{j-1})=1$, for each
$p\mid P_0$ the form $uq_j+vq_{j-1}$ vanishes modulo $p$ for exactly $p$ of the $p^2$ pairs $(u,v)$. By the Chinese remainder
theorem there are $1\le u,v\le P_0$ with $(k_0,2\Delta)=1$, where $k_0=uq_j+vq_{j-1}$. Put $h_0=up_j+vp_{j-1}$. Then
$|k_0\theta-h_0|\le2P_0/q_j$ and $k_0\le2P_0q_j$, so $|\theta-h_0/k_0|\le4P_0^2k_0^{-2}$. Write $h_0/k_0=h/k$ in lowest terms with
$k=k_0/g$; then $(k,2\Delta)=1$, and $(A+2)^{-1}k^{-2}\le|\theta-h/k|\le4P_0^2g^{-2}k^{-2}$ forces $g\le2P_0\sqrt{A+2}$.
\end{proof}

Apply Lemma~\ref{lem:hk} with $Y=C_0\sqrt{T/K}$, and put
\begin{equation}\label{eq:x0}
r=\frac hk,\qquad q=k\sqrt\Delta,\qquad\tau=\frac t{2\pi},\qquad x_0=\frac\tau r .
\end{equation}
Then $h\asymp k\asymp q\asymp\sqrt{T/K}$, and $r\asymp1$. For $|t-T|\le H$,
\[
x_0-x_*=\tau\Bigl(\frac1r-\sqrt\Delta\Bigr)+\frac{(t-T)\sqrt\Delta}{2\pi},\qquad
\Bigl|\frac1r-\sqrt\Delta\Bigr|=\frac{\sqrt\Delta}r\Bigl|r-\frac1{\sqrt\Delta}\Bigr|\le\frac{2\Delta c_4}{k^2}\le\frac{2\Delta c_4}{c_2^2C_0^2}\cdot\frac KT .
\]
Since $H\le T^{1/2-\delta+2\eps}\le K/(2\sqrt\Delta)$, we get $|x_0-x_*|\le K/2$ if $C_0$ is large enough. Hence all hypotheses of
Lemma~\ref{lem:statphase} hold, with $q\le T^{1/4}\le T^{1/2-\delta}$. In particular
\[
N_0=\frac{r^2q^2K^2}{x_*}=\frac{h^2\Delta K^2}{x_*}\asymp K .
\]
We fix $\eta(x)=\eta_0((x-x_*)/K)$ with $\eta_0$ as in \S\ref{sec:osc}; it satisfies the hypotheses of
Proposition~\ref{prop:zerodet}.

\subsection{The transformation formula}
Put
\begin{equation}\label{eq:varpi}
\varpi=\frac{\hb\Db}k-\frac1{2q^2r}=\frac{\hb\Db}k-\frac1{2hk\Delta},
\end{equation}
and, with $\phi(x)=\arsinh(x^{1/2})+(x+x^2)^{1/2}=2\sqrt x\,(1+\phi_1(x))$, where $\phi_1$ is analytic on $[0,1]$ and
$\phi_1(x)=x/6+O(x^2)$,
\begin{equation}\label{eq:Phi}
\Phi(u)=\frac t\pi\,\phi\Bigl(\frac{\pi u}{2q^2rt}\Bigr)=2\kappa\sqrt u\Bigl(1+\phi_1\Bigl(\frac{\pi^2\kappa^2u}{t^2}\Bigr)\Bigr),
\qquad\kappa=\Bigl(\frac t{2\pi q^2r}\Bigr)^{1/2}\asymp K^{1/2}.
\end{equation}

\begin{proposition}\label{prop:transf}
For $|t-T|\le H$,
\begin{multline*}
\Pc(t)=\e(\tau-\tau\log x_0)\sum_{\vs=\pm1}\sum_{w\in\Dc}\beta(w)\,\e\bigl(\varpi P(w)+\vs\Phi(P(w))\bigr)\mathcal A_\vs(P(w))\\
+O\bigl(T^{1/2}K^{-1}+KT^{-5/6}+(K/T)^{1/2}\bigr),
\end{multline*}
where $\mathcal A_\vs=\pi\theta_\vs q^{-1}\mathcal M_\vs$ satisfies $\mathcal A_\vs(u)=0$ for $u>5N_0$ and, for $u,N\ge1$,
\[
\mathcal A_\vs(u)\ll K^{1/4}u^{-1/4}T^{-1/2},\qquad\Var_{[N,2N]}\mathcal A_\vs\ll K^{1/4}N^{-1/4}T^{-1/2}.
\]
\end{proposition}
\begin{proof}
Write $\Pc(t)=\sum_nc(n)\e(-hn/k)f(n)$ with
\[
f(x)=\eta(x)x^{-1/2-it}\e\Bigl(\frac{hx}k\Bigr)=\eta(x)x^{-1/2}\e(g(x)),\qquad g(x)=rx-\tau\log x,
\]
and apply Corollary~\ref{cor:unified}. The integrals in \eqref{eq:unified} are $\mathcal J_\vs(P(w))$ of Lemma~\ref{lem:statphase}.
We treat the three parts of \eqref{eq:unified} in turn.

\emph{Main term.} By Corollary~\ref{cor:unified} and Lemma~\ref{lem:statphase}(iv), $\mathfrak m\int f\ll k^{-1}\ll(K/T)^{1/2}$.

\emph{$K_0$-terms.} By Lemma~\ref{lem:statphase}(iii), Corollary~\ref{cor:unified} and Lemma~\ref{lem:sector} (applied to $Q^*$),
\[
\mathfrak k(f)\ll\frac1k\sum_{n\ge1}\#\{w\in\Dc:P(w)=n\}\,KT^{-1/2}\e^{-c\sqrt{nT}/q}\ll\e^{-c'\sqrt K},
\]
since $\sqrt T/q\asymp\sqrt K$.

\emph{Dual terms.} By Lemma~\ref{lem:statphase}(i), each $\mathcal J_\vs(n)$ is the sum of
$\e(\tau-\tau\log x_0)\e(-n/(2q^2r)+\vs\Phi(n))\mathcal M_\vs(n)$ and $\mathcal E_\vs(n)$. In the first part,
$\e_k(\hb\Db n)\e(-n/(2q^2r))=\e(\varpi n)$, which gives the main term of the proposition. The second part contributes
\[
\ll\frac1q\sum_{\substack{w\in\Dc\\P(w)\le20N_0}}\frac{q^{1/2}}{P(w)^{1/4}T^{1/4}}\Bigl(\frac T{K^2}+T^{-1/3}\Bigr)
+\frac1q\sum_{\substack{w\in\Dc\\P(w)>20N_0}}T^{-A}\Bigl(\frac{P(w)}{N_0}\Bigr)^{-A}.
\]
Since $\#\{w\in\Dc:P(w)\le X\}\ll X$ for $X\ge1$, the first sum is $\ll q^{-1/2}T^{-1/4}N_0^{3/4}(T/K^2+T^{-1/3})$. With
$q^{-1/2}\asymp(K/T)^{1/4}$ and $N_0\asymp K$ this is $\ll KT^{-1/2}(T/K^2+T^{-1/3})=T^{1/2}K^{-1}+KT^{-5/6}$. The second sum is
negligible.

\emph{Amplitudes.} The properties of $\mathcal A_\vs$ follow from Lemma~\ref{lem:statphase}(ii), since $|\theta_\vs|=1$ and
$q^{-1}q^{1/2}T^{-1/4}\asymp K^{1/4}T^{-1/2}$.
\end{proof}

\begin{remark}
Proposition~\ref{prop:transf} is a smooth version of \cite[Lemma~3]{JS}, with a proof. The two sums over $j=1,2$ in
\cite[(3.6)]{JS} correspond to $\vs=\pm1$, and the cut-offs $n<n_j$ there to the condition $x_\vs(n)\in\supp\eta$ here
(Figure~\ref{fig:statpoints}). For indefinite forms the kernel $-Y_0$ gives $\theta_\vs=i\vs$ in place of $\theta_\vs=1$.
\end{remark}

\subsection{Reduction to exponential sums}
Let $\Dc_1,\dots,\Dc_J$ ($J\le2$) be the sets on which $\beta$ is constant: $\Dc_1=\Dc$ in the definite case, and
$\Dc_{1,2}=\mathcal C^*_\pm\cap\Z^2$ in the indefinite case. For $1\le N<N'\le2N$, $\vs=\pm1$ and $j\le J$ put
\begin{equation}\label{eq:SNdef}
S(N,N')=\sum_{\substack{w\in\Dc_j\\N\le P(w)\le N'}}\e\bigl(\varpi P(w)+\vs\Phi(P(w))\bigr).
\end{equation}

\begin{corollary}\label{cor:reduction}
Assume \eqref{eq:Krange}, and suppose that
\begin{equation}\label{eq:crucial2}
K^{1/4}N^{-1/4}T^{-1/2}|S(N,N')|\ll T^{-\eps}
\end{equation}
for all $|t-T|\le H$, $\vs$, $j$ and $1\le N\le5N_0$, $N<N'\le2N$. Then $\Pc(t)\ll(\log T)^{-2}$ for $|t-T|\le H$, and $Z$ has a zero
$\frac12+i\gamma$ with $|\gamma-T|\le H$.
\end{corollary}
\begin{proof}
Split the sum in Proposition~\ref{prop:transf} into the $J\le2$ sets $\Dc_j$ and into dyadic ranges $N\le P(w)<2N$, $N=2^i\le5N_0$.
For a function $\mathcal A$ of bounded variation on $[N,2N]$ and complex numbers $a_w$, partial summation gives
\[
\Bigl|\sum_{N\le P(w)<2N}a_w\mathcal A(P(w))\Bigr|\le\bigl(|\mathcal A(N)|+2\Var_{[N,2N]}\mathcal A\bigr)\max_{N\le N'<2N}\Bigl|\sum_{N\le P(w)\le N'}a_w\Bigr|.
\]
By Proposition~\ref{prop:transf} and \eqref{eq:crucial2}, each dyadic range contributes $\ll T^{-\eps}$, and there are $O(\log T)$
of them. The error terms in Proposition~\ref{prop:transf} are $\ll T^{-\delta}+T^{-7/30}+T^{-1/5}$ by \eqref{eq:Krange}.
Proposition~\ref{prop:zerodet} completes the proof.
\end{proof}

\section{Lemmas on exponential sums}\label{sec:lemmas}

Corollary~\ref{cor:reduction} reduces Theorems~\ref{thm:A}(a) and~\ref{thm:B} to a bound for the double exponential sum
$S(N,N')$. In this section we collect the tools for this bound: the B-process and van der Corput's estimate in one
variable, a uniform finiteness statement from o-minimality, Weyl's inequality in two variables, and a lemma
(Lemma~\ref{lem:2d}) which combines the one-variable estimates into a bound for sums over planar regions. In
\S\S\ref{sec:lemmas}--\ref{sec:expsum}, $\lambda$ and $\mu$ denote positive real parameters (not the linear forms of
\S\ref{subsec:automorphs}).

\begin{fact}[B-process; {\cite[Lemma~3.6]{GK}}, {\cite[Lemma~3]{BSS}}]\label{lem:B}
Let $f$ be real with four continuous derivatives on $[a,b]$, $b-a\le Y_0$, with $f''\asymp\lambda$ of constant sign,
$f'''\ll\lambda Y_0^{-1}$ and $f''''\ll\lambda Y_0^{-2}$. Let $f'(x_\nu)=\nu$. Then
\[
\sum_{a\le n\le b}\e(f(n))=\sum_{f'(a)\le\nu\le f'(b)}\frac{\e(f(x_\nu)-\nu x_\nu\pm\frac18)}{|f''(x_\nu)|^{1/2}}+O\bigl(\log(2+\lambda Y_0^2)+\lambda^{-1/2}\bigr),
\]
where $\nu$ runs over the integers between $f'(a)$ and $f'(b)$.
\end{fact}

\begin{fact}[van der Corput; {\cite[Thm.~2.2]{GK}}]\label{lem:vdc}
If $\mu\le|f''|\le\alpha\mu$ on $[a,b]$, then
\[
\sum_{a\le n\le b}\e(f(n))\ll\alpha(b-a)\mu^{1/2}+\mu^{-1/2}.
\]
\end{fact}

\begin{fact}[O-minimality; {\cite[Ch.~3]{vdD}}, {\cite{vdDM}}]\label{fact:omin}
Let $\R_{\mathrm{an}}$ be the expansion of the real field by all restricted analytic functions; it is o-minimal. If
$\mathcal F\subset\R^m\times\R^n$ is definable in $\R_{\mathrm{an}}$, there is $B$ such that for every $p\in\R^m$ the fibre
$\mathcal F_p=\{z:(p,z)\in\mathcal F\}$ has at most $B$ connected components. For $n=1$ these components are intervals, and,
by cell decomposition, the $j$-th of them is given by definable endpoint functions of $p$.
\end{fact}

\begin{lemma}[Two-dimensional Weyl inequality]\label{lem:weyl}
Let $a:\Z^2\to\C$ be supported on a finite set $\mathcal R$, and let $\mathcal M=\{m\in\Z^2:|m|_\infty<M\}$ with $M\ge1$. Then
\[
\Bigl|\sum_wa(w)\Bigr|^2\le\frac{\#(\mathcal R-\mathcal M)}{(2M-1)^2}\sum_{\ell\in\mathcal M-\mathcal M}\Bigl|\sum_wa(w+\ell)\overline{a(w)}\Bigr|.
\]
\end{lemma}
\begin{proof}
Write $(2M-1)^2\sum_wa(w)=\sum_{w\in\mathcal R-\mathcal M}\sum_{m\in\mathcal M}a(w+m)$, apply Cauchy--Schwarz, expand the square and
collect the terms with $m_1-m_2=\ell$. Each $\ell$ arises at most $(2M-1)^2$ times.
\end{proof}

The next lemma packages ``B-process in one variable, then van der Corput in the other''. In the language of
\cite[\S4]{BSS} these are the steps (4.5)--(4.25).

\begin{lemma}[Two-dimensional second-derivative estimate]\label{lem:2d}
Let $\mathcal A\subset\R^2$, $C_0\ge1$, $1\le X\le Y_0$ and $R\ge0$. Suppose that $\mathcal A$ meets at most $X$ of the lines
$\{x\}\times\R$, $x\in\Z$, and that for every $x\in\R$ the fibre $\mathcal A_x=\{y:(x,y)\in\mathcal A\}$ is a union of at most
$C_0$ intervals $I_{x,1}<\dots<I_{x,j(x)}$, each of length $\le Y_0$. Let $g$ be real and $C^4$ on a neighbourhood of
$\mathcal A$, and suppose:
\begin{enumerate}
\item[(i)] $\partial_y^2g\asymp\lambda$ with constant sign and $\|\nabla^2g\|\ll\lambda$ on $\mathcal A$; all third and fourth order
derivatives of $g$ are $\ll\lambda Y_0^{-1}$ and $\ll\lambda Y_0^{-2}$;
\item[(ii)] $\partial_yg(\mathcal A)$ is contained in an interval of length $R$;
\item[(iii)] $|\det\nabla^2g|/|\partial_y^2g|\asymp\mu$ on $\mathcal A$;
\item[(iv)] for each integer $\nu$ and $j\le C_0$, the set $X_{\nu,j}$ of those $x\in\R$ for which $\partial_yg(x,y)=\nu$ has a
solution $y=y_{\nu,j}(x)\in I_{x,j}$ is a union of $O(1)$ intervals, on each of which $y_{\nu,j}$ is continuous.
\end{enumerate}
Then
\[
\sum_{(x,y)\in\mathcal A\cap\Z^2}\e(g(x,y))\ll(R+1)\lambda^{-1/2}\bigl(X\mu^{1/2}+\mu^{-1/2}\bigr)+X\bigl(\log(2+\lambda Y_0^2)+\lambda^{-1/2}\bigr),
\]
where the implied constant depends on $C_0$ and on the constants in (i)--(iv).
\end{lemma}
\begin{proof}
Since $\partial_y^2g$ has constant sign, $\partial_yg(x,\cdot)$ is strictly monotone on each $I_{x,j}$, so $y_{\nu,j}(x)$ is unique.
For each $(x,j)$, Fact~\ref{lem:B} gives
\[
\sum_{y\in I_{x,j}}\e(g(x,y))=\sum_{\nu\in\partial_yg(x,I_{x,j})}\frac{\e(G_{\nu,j}(x)\pm\frac18)}{|\partial_y^2g(x,y_{\nu,j}(x))|^{1/2}}+O\bigl(\log(2+\lambda Y_0^2)+\lambda^{-1/2}\bigr),
\]
with $G_{\nu,j}(x)=g(x,y_{\nu,j}(x))-\nu y_{\nu,j}(x)$. Summed over the at most $C_0X$ pairs $(x,j)$, the error terms give the
last term. In the main terms we interchange the summations:
\[
\sum_x\sum_{j\le j(x)}\sum_{\nu\in\partial_yg(x,I_{x,j})}=\sum_{j\le C_0}\sum_\nu\sum_{x\in X_{\nu,j}} .
\]
\emph{Counting the frequencies.} Every $\nu$ that occurs lies in $\partial_yg(\mathcal A)$, so by (ii) at most $R+1$ integers $\nu$
occur in total. (A bound for the variation of $\partial_yg$ along each fibre would not suffice here.)

\emph{The sum over $x$.} Fix $j,\nu$ and one of the intervals $\mathcal I$ of (iv); we have to bound the sum over
$x\in\mathcal I\cap\Z$. Since every such $x$ is one of the at most $X$ integers with $\mathcal A_x\ne\emptyset$, we may replace
$\mathcal I$ by the convex hull of $\mathcal I\cap\Z$, which has length $<X$. By (i) and the implicit function theorem,
$y_{\nu,j}$ is $C^3$ on $\mathcal I$ with $y'_{\nu,j}=-\partial_x\partial_yg/\partial_y^2g\ll1$. Moreover
\[
G'_{\nu,j}(x)=\partial_xg(x,y_{\nu,j}(x)),\qquad G''_{\nu,j}=\frac{\det\nabla^2g}{\partial_y^2g},
\]
so $|G''_{\nu,j}|\asymp\mu$ by (iii). The amplitude $|\partial_y^2g(x,y_{\nu,j}(x))|^{-1/2}$ is $\asymp\lambda^{-1/2}$, and its derivative is
$\ll\lambda^{-1/2}Y_0^{-1}$, so its total variation on $\mathcal I$ ($|\mathcal I|<X\le Y_0$) is $\ll\lambda^{-1/2}$. Partial summation and
Fact~\ref{lem:vdc} bound the sum over $x\in\mathcal I$ by $\ll\lambda^{-1/2}(X\mu^{1/2}+\mu^{-1/2})$. Multiplying by the $O(1)$ intervals,
the $\le C_0$ values of $j$ and the $\le R+1$ values of $\nu$ proves the lemma.
\end{proof}

\section{The double exponential sum}\label{sec:expsum}

We now prove the key estimate, Theorem~\ref{thm:SN}. The proof follows the plan of \S\ref{subsec:step3}: Weyl differencing in
two variables (\S\ref{subsec:proofSN}, Step~1), an analysis of the Hessian of the differenced phase
(\S\ref{subsec:diffphase}), a decomposition of the summation region according to the angle between $w$ and the shift
(\S\ref{subsec:decomp}), and Lemma~\ref{lem:2d} on each piece.

\begin{setting}\label{setting}
Corollary~\ref{cor:reduction} leads to the following situation, which we fix for this section.
\begin{itemize}
\item $P$ is an integral binary quadratic form with discriminant $d\ne0$, $B(v,w)=P(v+w)-P(v)-P(w)$, and $\mathcal C$ is
one of the following:
\begin{itemize}
\item[--] \emph{definite case:} $P$ positive definite and $\mathcal C=\mathcal C^\sharp=\R^2\setminus\{0\}$;
\item[--] \emph{indefinite case:} $d>0$ not a square, and $\mathcal C$ a fundamental cone as in Lemma~\ref{lem:sector}, so
that $\overline{\mathcal C}\setminus\{0\}\subset\{P>0\}$; $\mathcal C^\sharp\supset\overline{\mathcal C}\setminus\{0\}$ is a slightly
larger open convex cone with $\overline{\mathcal C^\sharp}\setminus\{0\}\subset\{P>0\}$.
\end{itemize}
\item $T,K$ satisfy \eqref{eq:Krange}, $t\asymp T$, $\kappa\asymp K^{1/2}$, $\Phi$ is given by \eqref{eq:Phi}, $\varpi\in\R$
and $\vs\in\{\pm1\}$.
\item $c_0K^{1/2}\le N\le C_1K$, $N<N'\le2N$, $L=N^{1/2}$, and
\[
\mathcal R=\{w\in\Z^2\cap\mathcal C:N\le P(w)\le N'\},\qquad f(w)=\varpi P(w)+\vs\Phi(P(w)).
\]
\end{itemize}
In both cases $P(w)\asymp|w|^2$ on $\mathcal C^\sharp$, so $|w|\asymp L$ on $\mathcal R$, and $L\le C\kappa$. All implied constants in
this section depend only on $P$, $\mathcal C$, $\mathcal C^\sharp$, $c_0$, $C_1$, $\delta$ and the constants implicit in $t\asymp T$,
$\kappa\asymp K^{1/2}$.
\end{setting}

In the application, $P=Q^*$ (definite case) or $P=\pm Q^*$ with $\mathcal C=\mathcal C^*_\pm$ (indefinite case), and
$d=\operatorname{disc}Q^*=b^2-4ac$.

\begin{theorem}\label{thm:SN}
In Setting~\ref{setting},
\[
S(N)=\sum_{w\in\mathcal R}\e(f(w))\ll N^{2/3}K^{1/6}\log T .
\]
\end{theorem}

\subsection{Derivatives of $\sqrt P$}
Let $J=\bigl(\begin{smallmatrix}0&-1\\1&0\end{smallmatrix}\bigr)$, so that $\det(v,w)=(Jv)\cdot w$. For $w,e\in\R^2$,
\begin{equation}\label{eq:Bdisc}
B(w,e)^2-4P(w)P(e)=d\det(w,e)^2 .
\end{equation}
On $\{P>0\}$ we therefore have
\begin{equation}\label{eq:der}
\nabla^2\sqrt P(w)=-\frac d4P(w)^{-3/2}(Jw)(Jw)^T,\qquad\partial_e^3\sqrt P(w)=\frac38\,d\det(w,e)^2B(w,e)P(w)^{-5/2}.
\end{equation}
Indeed, let $p(s)=P(w+se)=P(w)+sB(w,e)+s^2P(e)$ and $q=p^{1/2}$. Then
\[
q''=\frac{4P(e)p-p'^2}{4p^{3/2}},
\]
whose numerator has derivative $4P(e)p'-2p'p''=0$; at $s=0$ the numerator is $-d\det(w,e)^2$ by \eqref{eq:Bdisc}. Since this
holds for every $e$, and $(Jw\cdot e)^2=\det(w,e)^2$, the Hessian is $-\frac d4P^{-3/2}(Jw)(Jw)^T$. In particular the Hessian of
$\sqrt P$ has rank one and kernel $\R w$; this reflects the homogeneity of degree one of $\sqrt P$.

\begin{lemma}\label{lem:dir}
Let $E=\{(1,0),(0,1),(1,1)\}$. There is $c_E>0$ such that every real symmetric $2\times2$ matrix $H$ satisfies
$\max_{e\in E}|e^THe|\ge c_E\|H\|$.
\end{lemma}
\begin{proof}
A binary quadratic form vanishing at three pairwise independent vectors is zero; now use compactness of the unit sphere
of symmetric matrices.
\end{proof}

\subsection{The differenced phase}\label{subsec:diffphase}
For $\ell\in\R^2$ put
\[
g_\ell(w)=\Phi(P(w+\ell))-\Phi(P(w)),\qquad\lambda=\lambda(\ell)=\frac{\kappa|\ell|}{L^2}.
\]
Let $\psi\in[0,\frac\pi2]$ be the angle between the lines $\R w$ and $\R\ell$, and let $\Hs_\ell(w)$ be the Hessian of $g_\ell$.

\begin{lemma}\label{lem:hess}
There is $C>0$ such that for $0<|\ell|\le L/C$ and all $w\in\mathcal C^\sharp$ with $c'L\le|w|\le C'L$ the following holds
(with $c',C'$ fixed): $w+\ell\in\mathcal C^\sharp$, and
\begin{enumerate}
\item[(a)] $\|\Hs_\ell(w)\|\asymp\lambda$, and all third and fourth order derivatives of $g_\ell$ are $\ll\lambda L^{-1}$ and
$\ll\lambda L^{-2}$, respectively;
\item[(b)] if $\sin^2\psi\ge T^{-3/5}$, then $|\det\Hs_\ell(w)|\asymp(\lambda\sin\psi)^2$;
\item[(c)] in the orthonormal basis $(\nu,\hat w)$, $\hat w=w/|w|$, $\nu=J\hat w$, write
$\Hs_\ell=\bigl(\begin{smallmatrix}a_{11}&a_{12}\\a_{12}&a_{22}\end{smallmatrix}\bigr)$. Then
\[
|a_{12}|\ll\lambda\sin\psi+\lambda T^{-4/5},\qquad|a_{22}|\ll\lambda\sin^2\psi+\lambda T^{-4/5}.
\]
\end{enumerate}
\end{lemma}
\begin{proof}
The first claim holds for $C$ large, since $\mathcal C^\sharp$ contains a fixed angular neighbourhood of
$\overline{\mathcal C}\setminus\{0\}$ (definite case: trivially).

\emph{The model phase.} By \eqref{eq:Phi}, $\Phi(P)=2\kappa\sqrt P+F_1$ with $F_1=2\kappa\sqrt P\,\phi_1(\pi^2\kappa^2P/t^2)$. Put
\[
\varrho=\frac{\pi^2\kappa^2L^2}{t^2}\asymp\frac{KN}{T^2}\le C_1\frac{K^2}{T^2}\ll T^{-4/5}.
\]
Every derivative of order $j$ of $F_1$ is $\ll\kappa L^{1-j}\varrho$ on $|w|\asymp L$, $w\in\mathcal C^\sharp$. Hence the second, third
and fourth derivatives of $F_1(w+\ell)-F_1(w)$ are $\ll\lambda\varrho$, $\lambda\varrho L^{-1}$ and $\lambda\varrho L^{-2}$. It therefore
suffices to prove (a)--(c) for
\[
g^0_\ell(w)=2\kappa\bigl(\sqrt{P(w+\ell)}-\sqrt{P(w)}\bigr)
\]
and its Hessian $\Hs^0_\ell$, with the terms $\lambda T^{-4/5}$ omitted in (c). Indeed, the perturbation $E=\Hs_\ell-\Hs^0_\ell$,
$\|E\|\ll\lambda T^{-4/5}$, changes the determinant by $O(\lambda^2T^{-4/5})$, which is admissible in (b) because
$\sin^2\psi\ge T^{-3/5}$.

\emph{Determinant.} By \eqref{eq:der},
\[
\Hs^0_\ell(w)=-\frac{\kappa d}2\bigl(\alpha\,aa^T-\beta\,bb^T\bigr),\qquad a=J(w+\ell),\ b=Jw,\ \alpha=P(w+\ell)^{-3/2},\ \beta=P(w)^{-3/2}.
\]
For $2\times2$ matrices $\det(\alpha aa^T-\beta bb^T)=-\alpha\beta\det[a,b]^2$, and $\det[a,b]=\det[\ell,w]$, $|\det[\ell,w]|=|\ell||w|\sin\psi$.
This gives $|\det\Hs^0_\ell|\asymp(\lambda\sin\psi)^2$ \emph{exactly} up to the constants in $\alpha,\beta\asymp L^{-3}$.

\emph{Norm.} The identity
\[
\alpha aa^T-\beta bb^T=(\alpha-\beta)bb^T+\alpha\bigl(b(J\ell)^T+(J\ell)b^T+(J\ell)(J\ell)^T\bigr),\qquad|\alpha-\beta|\ll|\ell|L^{-4},
\]
gives the upper bound. If $\sin\psi\ge c$, the lower bound follows from the determinant. If $\sin\psi<c$, write $\ell=sw+\ell'$ with
$|\ell'|\le c|\ell|$, so $a=(1+s)b+J\ell'$ and $\alpha=(1+s)^{-3}\beta(1+O(c|\ell|/L))$. The coefficient of $bb^T$ is then
$\beta((1+s)^{-1}-1)+O(c|\ell|L^{-4})$, and $|s|\asymp|\ell|/L$; the remaining terms have norm $\ll c|\ell|L^{-2}$. For $c$ small this
gives $\|\Hs^0_\ell\|\gg\lambda$.

\emph{Entries.} Since $b\parallel\nu$, we have $b\cdot\hat w=0$ and $a\cdot\hat w=(J\ell)\cdot\hat w=\pm|\ell|\sin\psi$, while $|a\cdot\nu|\ll L$. Hence
\[
a_{12}=-\frac{\kappa d}2\alpha(a\cdot\nu)(a\cdot\hat w)\ll\kappa L^{-3}L|\ell|\sin\psi=\lambda\sin\psi,\qquad
a_{22}=-\frac{\kappa d}2\alpha(a\cdot\hat w)^2\ll\lambda\sin^2\psi .
\]
\emph{Higher derivatives.} $\sqrt P$ is homogeneous of degree one, so its derivatives of order $j\ge1$ are $\ll L^{1-j}$ on
$|w|\asymp L$, $w\in\mathcal C^\sharp$; apply the mean value theorem to $g^0_\ell$.
\end{proof}

\subsection{Decomposition relative to a shift}\label{subsec:decomp}
Fix $\ell\in\Z^2$ with $0<|\ell|\le2L^{2/3}$ (so $|\ell|\le L/C$ for $N$ large), and put $\Psi_0=(\kappa|\ell|)^{-1/2}$. Since
$L^{2/3}\le C\kappa^{2/3}\ll K^{1/3}$,
\[
\Psi_0^2\gg(K^{1/2}K^{1/3})^{-1}\ge T^{-1/2}.
\]
Let $\mathcal R_\ell=\mathcal R\cap(\mathcal R-\ell)$. For dyadic $\Psi=2^i\Psi_0\le1$ put
\[
\mathcal A(\Psi)=\{w\in\mathcal R_\ell:\Psi<\psi(w)\le2\Psi\}.
\]
On $\mathcal A(\Psi)$ we have $\sin\psi\asymp\Psi$ and $\sin^2\psi\gg\Psi_0^2\ge T^{-3/5}$, so Lemma~\ref{lem:hess}(b) applies.

\begin{lemma}\label{lem:decomp}
Each $\mathcal A(\Psi)$ is the disjoint union of $O(1)$ sets $\mathcal A$ (the bound being uniform in $N,N',\ell,\Psi,\varpi,t,\kappa$),
each of which is the intersection of $\Z^2$ with a set $\mathcal A^\R\subset\R^2$, such that for some $e\in E$:
\begin{enumerate}
\item[(a)] $|\partial_e^2g_\ell|\asymp\lambda$ with constant sign on $\mathcal A^\R$;
\item[(b)] $|\det\Hs_\ell|\asymp(\lambda\Psi)^2$ on $\mathcal A^\R$;
\item[(c)] $\sup_{\mathcal A^\R}\partial_eg_\ell-\inf_{\mathcal A^\R}\partial_eg_\ell\ll\lambda L\Psi$;
\item[(d)] every line parallel to $e$ meets $\mathcal A^\R$ in $O(1)$ intervals, each of length $\ll L\Psi$;
\item[(e)] $\mathcal A^\R\subset\{|w|\le C'L\}$.
\end{enumerate}
\end{lemma}
\begin{proof}
Let $\mathcal R^\R$ be the real region defined by the inequalities defining $\mathcal R$ (without $w\in\Z^2$), and
\[
\mathcal R^\R_\ell=\mathcal R^\R\cap(\mathcal R^\R-\ell),\qquad\mathcal A^\R(\Psi)=\{w\in\mathcal R^\R_\ell:\Psi<\psi(w)\le2\Psi\}.
\]
Fix a small constant $c''>0$, to be chosen.

\emph{Uniformity.} After the scaling $w=Lu$, the sets and functions below form families definable in $\R_{\mathrm{an}}$ with
finitely many real parameters: $\vartheta=|\ell|/L$, $\hat\ell=\ell/|\ell|$, $\varrho$, $N/L^2$, $N'/L^2$ and $\Psi$. Indeed
\begin{equation}\label{eq:scaled}
\frac{g_\ell(Lu)}{\kappa|\ell|}=\frac2{\vartheta}\Bigl(\sqrt{P(u+\vartheta\hat\ell)}\,\bigl(1+\phi_1(\varrho P(u+\vartheta\hat\ell))\bigr)-\sqrt{P(u)}\,\bigl(1+\phi_1(\varrho P(u))\bigr)\Bigr),
\end{equation}
where $\phi_1$ is only evaluated on $[0,1]$, and $\sqrt\cdot$ only on $[c,\infty)$. By Fact~\ref{fact:omin}, every set constructed
below from $g_\ell$, its derivatives and $\mathcal A^\R(\Psi)$ by finitely many Boolean operations and inequalities has a
number of connected components bounded uniformly in the parameters.

\emph{Small $\Psi$ ($\Psi\le c''$).} By Lemma~\ref{lem:hess}(a),(c), $|a_{12}|\ll\lambda\Psi$ and $|a_{22}|\ll\lambda\Psi^2$ (as
$T^{-4/5}\le\Psi_0^2$), so $|a_{11}|\ge\|\Hs_\ell\|-O(\lambda\Psi)\gg\lambda$ if $c''$ is small. For each unit vector $\hat w$ at least two
elements $e\in E$ satisfy $|e\cdot J\hat w|\ge c_2|e|$; by compactness we can cover the unit circle by finitely many closed
arcs $I_1,\dots,I_r$ and choose $e_i\in E$ with $|e_i\cdot J\hat w|\ge c_2|e_i|$ for $\hat w\in I_i$. Take for the sets $\mathcal A^\R$
the connected components of $\{w\in\mathcal A^\R(\Psi):\hat w\in I_i\}\setminus\bigcup_{i'<i}\{\hat w\in I_{i'}\}$, with $e=e_i$. Then
\[
e^T\Hs_\ell e=a_{11}(e\cdot\nu)^2+2a_{12}(e\cdot\nu)(e\cdot\hat w)+a_{22}(e\cdot\hat w)^2=a_{11}(e\cdot\nu)^2+O(\lambda\Psi),
\]
so $|\partial_e^2g_\ell|\asymp\lambda$; being continuous and non-zero on a connected set, it has constant sign. This gives (a), and
(b) is Lemma~\ref{lem:hess}(b). Each $\mathcal A^\R$ lies in one of the four sectors $\{\Psi<\psi\le2\Psi\}$ adjacent to the
line $\R\ell$, of opening $\Psi$, and in $\{|w|\le C'L\}$; in the indefinite case this sector lies in $\mathcal C^\sharp$ for $c''$
small, because it meets $\mathcal C$ and has opening $\Psi\le c''$. A line parallel to $e$ makes an angle $\gg1$ with every
direction in that sector, so it meets the sector in a segment of length $\ll L\Psi$; the number of intervals is $O(1)$ by
Fact~\ref{fact:omin}. This gives (d) and (e). For (c), let $\widetilde{\mathcal A}$ be the part of that sector with
$c'L\le|w|\le C'L$; it contains $\mathcal A^\R$ and lies in $\mathcal C^\sharp$. Two points of $\widetilde{\mathcal A}$ are joined
inside $\widetilde{\mathcal A}$ by a radial segment of length $\ll L$, on which
\[
\partial_{\hat w}\partial_eg_\ell=(e\cdot\nu)a_{12}+(e\cdot\hat w)a_{22}\ll\lambda\Psi,
\]
followed by a circular arc of length $\ll L\Psi$, on which $|\nabla\partial_eg_\ell|\ll\lambda$. So $\partial_eg_\ell$ varies by
$\ll\lambda L\Psi$ on $\widetilde{\mathcal A}$.

\emph{Large $\Psi$ ($c''<\Psi\le1$).} There are $O(1)$ such $\Psi$. For $e\in E$ let
\[
U_e=\{w\in\mathcal A^\R(\Psi):|e^T\Hs_\ell e|\ge c_E\|\Hs_\ell\|\};
\]
these cover $\mathcal A^\R(\Psi)$ by Lemma~\ref{lem:dir}. Take for the sets $\mathcal A^\R$ the connected components of $U_{e_1}$,
$U_{e_2}\setminus U_{e_1}$ and $U_{e_3}\setminus(U_{e_1}\cup U_{e_2})$, with the corresponding $e$. Then (a) holds by
Lemma~\ref{lem:hess}(a) and connectedness, (b) by Lemma~\ref{lem:hess}(b), (d) and (e) by Fact~\ref{fact:omin} and
$\mathcal A^\R\subset\{|w|\le C'L\}$ (lengths $\ll L\asymp L\Psi$). For (c), the region $\{w\in\mathcal C^\sharp:c'L\le|w|\le C'L\}$ is
connected by paths of length $\ll L$ (radial segments and arcs), and $|\nabla\partial_eg_\ell|\ll\lambda$ on it.
\end{proof}

\subsection{Proof of Theorem~\ref{thm:SN}}\label{subsec:proofSN}
\emph{Step 1: Weyl differencing.} Let $1\le M\le L^{2/3}/2$, so that the shifts below satisfy $|\ell|\le3M\le2L^{2/3}$. Since
$\mathcal R\subset\{|w|\le C'L\}$, we have $\#(\mathcal R-\mathcal M)\ll L^2$. Lemma~\ref{lem:weyl} with $a(w)=\mathbf 1_{\mathcal R}(w)\e(f(w))$ gives
\begin{equation}\label{eq:afterweyl}
|S(N)|^2\ll\frac{L^2}{M^2}\Bigl(L^2+\sum_{\substack{\ell\in\Z^2\\0<|\ell|\le3M}}|S_\ell|\Bigr),\qquad
S_\ell=\sum_{w\in\mathcal R_\ell}\e\bigl(f(w+\ell)-f(w)\bigr).
\end{equation}
Here
\[
f(w+\ell)-f(w)=\varpi B(w,\ell)+\varpi P(\ell)+\vs g_\ell(w).
\]
The contribution of $\varpi P$ is \emph{linear} in $w$; in particular the rational part of $\varpi$ no longer influences any
second derivative.

\emph{Step 2: the thin cone.} The set $\{w\in\mathcal R_\ell:\psi(w)\le\Psi_0\}$ lies in the union of two convex sectors of
opening $\ll\Psi_0$ and radius $\ll L$. Each has area $\ll L^2\Psi_0$ and perimeter $\ll L$, so it contains $\ll L^2\Psi_0+L$
lattice points.

\emph{Step 3: the pieces.} Let $\mathcal A=\Z^2\cap\mathcal A^\R$ be one of the sets of Lemma~\ref{lem:decomp}. Choose $G\in SL_2(\Z)$
with $Ge_2=e$ from a fixed finite set (one for each $e\in E$), and substitute $w=Gz$, $z\in\Z^2$. The linear part of the phase
stays linear in $z$; it shifts $\partial_{z_2}$ of the phase by a constant and does not affect the hypotheses below. Since
$\partial_{z_2}=\partial_e$ and $\det G=1$, Lemmas~\ref{lem:decomp} and \ref{lem:hess} give the hypotheses of Lemma~\ref{lem:2d} for
$g(z)=\vs g_\ell(Gz)$ with
\[
\lambda\asymp\kappa|\ell|L^{-2},\quad\mu\asymp\lambda\Psi^2,\quad R\ll\lambda L\Psi,\quad Y_0\asymp L,\quad X\ll L,
\]
the constants depending on the finite set of $G$. Hypothesis (iv) holds by Fact~\ref{fact:omin}: by \eqref{eq:scaled}, the sets
$X_{\nu,j}$ and functions $y_{\nu,j}$ are fibres of families definable in $\R_{\mathrm{an}}$, with the additional parameter
$\nu L/(\kappa|\ell|)$ (shifted by the constant from the linear part). Lemma~\ref{lem:2d} gives
\begin{align*}
\sum_{w\in\mathcal A}\e(\cdots)&\ll(\lambda L\Psi+1)\lambda^{-1/2}\bigl(L\lambda^{1/2}\Psi+\lambda^{-1/2}\Psi^{-1}\bigr)+L\bigl(\log T+\lambda^{-1/2}\bigr)\\
&\ll\kappa|\ell|\Psi^2+L\log T+\frac{L^2}{\kappa|\ell|\Psi}+\frac{L^2}{(\kappa|\ell|)^{1/2}} ,
\end{align*}
using $\lambda L^2=\kappa|\ell|$ and $\lambda L\Psi\cdot\lambda^{-1}\Psi^{-1}=L$. Summing over the $O(1)$ sets $\mathcal A$ in each $\mathcal A(\Psi)$
and over the $O(\log T)$ values of $\Psi\ge\Psi_0$ (note $L^2/(\kappa|\ell|\Psi)\le L^2/(\kappa|\ell|)^{1/2}$), and adding Step~2,
\begin{equation}\label{eq:Sl}
|S_\ell|\ll\Bigl(\kappa|\ell|+L+\frac{L^2}{(\kappa|\ell|)^{1/2}}\Bigr)\log^2T .
\end{equation}

\emph{Step 4: choice of $M$.} The number of $\ell\in\Z^2$ with $|\ell|\le3M$ is $\ll M^2$, and $\sum_\ell|\ell|\ll M^3$,
$\sum_\ell|\ell|^{-1/2}\ll M^{3/2}$. Inserting \eqref{eq:Sl} into \eqref{eq:afterweyl},
\[
|S(N)|^2\ll\Bigl(\frac{L^4}{M^2}+L^2\kappa M+L^3+\frac{L^4}{(\kappa M)^{1/2}}\Bigr)\log^2T .
\]
Choose $M=\lfloor c_3(L^2/\kappa)^{1/3}\rfloor$ with $c_3>0$ small. Since $N\ge c_0K^{1/2}$ and $\kappa\asymp K^{1/2}$, a suitable
choice of $c_0$ gives $M\ge1$; and $M\le L^{2/3}/2$ because $\kappa\gg1$. The first two terms are then $\asymp L^{8/3}\kappa^{2/3}$.
The third is smaller because $L\le\kappa^2$. The fourth is $\asymp L^{11/3}\kappa^{-1/3}$, which is $\ll L/\kappa\ll1$ times the first.
Hence
\[
|S(N)|\ll L^{4/3}\kappa^{1/3}\log T\asymp N^{2/3}K^{1/6}\log T .\qquad\qed
\]

\section{Proof of Theorems~\ref{thm:A}(a) and~\ref{thm:B}}\label{sec:proof}

Let $K=T^{3/5-2\eps}$ and $H=T^{1+2\eps}/K=T^{2/5+4\eps}$; then \eqref{eq:Krange} holds with $\delta=\frac1{20}$ for $\eps\le\frac1{80}$.
Choose $h/k$ as in \S\ref{subsec:twist}. By Corollary~\ref{cor:reduction} it suffices to prove \eqref{eq:crucial2} for
$1\le N\le5N_0$. The sum $S(N,N')$ of \eqref{eq:SNdef} is the sum $S(N)$ of Setting~\ref{setting}, with $P=Q^*$,
$\mathcal C=\R^2\setminus\{0\}$ (definite case) or $P=\pm Q^*$, $\mathcal C=\mathcal C^*_\pm$ (indefinite case), and with $C_1K\ge5N_0$.

For $c_0K^{1/2}\le N\le5N_0$, Theorem~\ref{thm:SN} gives
\[
K^{1/4}N^{-1/4}T^{-1/2}|S(N,N')|\ll K^{5/12}N^{5/12}T^{-1/2}\log T\ll K^{5/6}T^{-1/2}\log T=T^{-5\eps/3}\log T .
\]
For $N<c_0K^{1/2}$ the trivial bound $|S(N,N')|\ll N$ gives $K^{1/4}N^{3/4}T^{-1/2}\ll K^{5/8}T^{-1/2}\le T^{-1/8}$. Hence
\eqref{eq:crucial2} holds (with $\eps$ replaced by a smaller positive constant), and Corollary~\ref{cor:reduction} gives a zero
$\frac12+i\gamma$ of $Z$ with $|\gamma-T|\le T^{2/5+4\eps}$. Since $\eps$ is arbitrary, this proves Theorem~\ref{thm:A}(a) and Theorem~\ref{thm:B}.\qed

\section{Positive definite forms with square $\Delta$}\label{sec:square}

In this section we prove Theorem~\ref{thm:A}(b), and we explain why the method fails for $\Delta=4m^2$ with $m$ even.

Let $Q=[a,b,c]$ be primitive and positive definite with $\Delta=4ac-b^2=f^2$. Since $\Delta\equiv-b^2\equiv0,3\pmod4$, $b$ and $f$ are
even; write $b=2b'$ and $f=2m$. These are the forms of discriminant $-4m^2$, attached to orders in $\Z[i]$. Jutila and
Srinivas exclude this case (\cite[\S1]{JS}: ``other cases like $\Delta=4$ \dots\ are easier or well-known''). Here
$1/\sqrt\Delta=1/f$ is rational, and we take $h/k=1/f$ itself. Then $q=k\sqrt\Delta=f^2$ is bounded.

\begin{lemma}\label{lem:Gsq}
Let $(h,f)=1$. If $m$ is odd, then $\GQ(f,h)=0$. If $m$ is even, then $|\GQ(f,h)|=f\sqrt{2f}$.
\end{lemma}
\begin{proof}
Put $A'=\bigl(\begin{smallmatrix}a&b'\\b'&c\end{smallmatrix}\bigr)$, so that $A=2A'$, $\det A'=m^2$ and $Q(v)=v^TA'v$. Then
$\mathcal R_f=\{v\bmod2m:A'v\equiv0\pmod m\}$, and by Lemma~\ref{lem:gauss}(b) with $w=0$, $\GQ(f,h)\ne0$ if and only if
$Q(v)\equiv0\pmod{2m}$ for all $v\in\mathcal R_f$, in which case $|\GQ(f,h)|=f|\mathcal R_f|^{1/2}$. For $v\in\mathcal R_f$ we always have
$Q(v)=v\cdot A'v\equiv0\pmod m$. Since $Q$ is primitive and $b$ is even, $a$ or $c$ is odd.

\emph{$m$ odd.} The vectors $(m,0)$ and $(0,m)$ lie in $\mathcal R_f$, and $Q$ takes the values $m^2a$, $m^2c$ there, one of which
is odd. So $\GQ(f,h)=0$.

\emph{$m$ even.} Let $2^j\parallel m$, $j\ge1$, and suppose $a$ odd (otherwise exchange the variables). For $v\in\mathcal R_f$ we
have $av_1+b'v_2\equiv0\pmod{2^j}$, so
\[
aQ(v)=(av_1+b'v_2)^2+m^2v_2^2\equiv0\pmod{2^{2j}} .
\]
Hence $Q(v)\equiv0\pmod{2^{j+1}}$, and together with $Q(v)\equiv0\pmod m$ this gives $Q(v)\equiv0\pmod{2m}$. Finally, the condition
$A'v\equiv0\pmod m$ depends on $v\bmod m$ only, and the Smith form of $A'$ is $\mathrm{diag}(1,m^2)$ since $\gcd(a,b',c)=1$; so
$|\mathcal R_f|=4m=2f$.
\end{proof}

\begin{theorem}[Theorem~\ref{thm:A}(b)]\label{thm:square}
Let $Q$ be primitive and positive definite with $\Delta=4m^2$, $m$ odd. Then for every $\eps>0$ and $T\ge T_0(\eps,Q)$, $\zeta_Q$ has a
zero $\frac12+i\gamma$ with $|\gamma-T|\le T^{1/3+\eps}$.
\end{theorem}
\begin{proof}
Let $K=T^{2/3-2\eps}$ and $H=T^{1+2\eps}/K=T^{1/3+4\eps}$, and let $\eta$ be as in \S\ref{subsec:twist}. Write
$\Pc(t)=\sum_nr_Q(n)\e(-n/f)F(n)$ with $F(x)=\eta(x)x^{-1/2}\e(g(x))$, $g(x)=x/f-\tau\log x$, and apply
Proposition~\ref{prop:Vdef} with $k=f$ and $h=-1$:
\[
\Pc(t)=\frac{2\pi\GQ(f,-1)}{f^3}\int_0^\infty F+\frac\pi{f^3}\sum_{w\ne0}\GQ(f,-1,w)\sum_{\vs=\pm1}\mathcal J_\vs(Q^*(w)),
\]
with $\mathcal J_\vs$ as in Lemma~\ref{lem:statphase} for $r=1/f$, $q=f^2$. The hypotheses of that lemma hold: $x_0=tf/2\pi$ and
$x_*=Tf/2\pi$ differ by at most $fH/2\pi\le K/2$. The first term vanishes by Lemma~\ref{lem:Gsq}. By Lemma~\ref{lem:gauss}(b),
$|\GQ(f,-1,w)|\le f\sqrt{2f}\ll1$. Here $N_0=f^2K^2/x_*\asymp K^2/T$, and Lemma~\ref{lem:statphase}(i),(ii) gives
\[
\Pc(t)\ll\sum_{\substack{w\ne0\\Q^*(w)\le20N_0}}Q^*(w)^{-1/4}T^{-1/4}+T^{-1}\ll T^{-1/4}N_0^{3/4}\ll K^{3/2}T^{-1}=T^{-3\eps}.
\]
Proposition~\ref{prop:zerodet} completes the proof.
\end{proof}

For even $m$ the same choice $h/k=1/f$ fails, and not for lack of a good bound for the dual sum: the main term of the
transformation contributes a non-oscillating constant to the Hardy function.

\begin{proposition}\label{prop:square}
Let $Q$ be primitive and positive definite with $\Delta=4m^2$, $m$ even, and let
$W(t)=\e^{\pi t/2}(\sqrt\Delta/2\pi)^{1/2+it}\Ga(\frac12+it)\zeta_Q(\frac12+it)$ be the Hardy function of \cite[\S2]{JS}. With
$K,H,\eta$ as in Theorem~\ref{thm:square} (or any $T^{1/2+\delta}\le K\le T^{1-\delta}$ with $fH\le\pi K$), the main term
$\frac{2\pi\GQ(f,-1)}{f^3}\int F$ of the transformation of $\Pc(t)$, multiplied by the gamma factor of $W$, equals
\[
c_Q\Bigl(1+O\Bigl(\frac T{K^2}+T^{-1/3}\Bigr)\Bigr),\qquad c_Q=\frac{2\pi}{f^2}\,\e^{i\pi/4}\GQ(f,-1),\qquad|c_Q|=2\pi\sqrt{\frac2f}.
\]
Consequently it contributes $\sqrt\pi\,c_QH_0(1+o(1))$ to the integral $I$ of \cite[\S2]{JS}.
\end{proposition}
\begin{proof}
By Lemma~\ref{lem:statphase}(iv) with $r=1/f$ and $\eta(x_0)=1$,
\[
\int_0^\infty F=\sqrt f\,\e^{i\pi/4}\e^{i(t-t\log(tf/2\pi))}\Bigl(1+O\Bigl(\frac T{K^2}+T^{-1/3}\Bigr)\Bigr).
\]
By Stirling's formula, $\Ga(\frac12+it)=\sqrt{2\pi}\,\e^{-\pi t/2}\e^{i(t\log t-t)}(1+O(t^{-1}))$, so
\[
\e^{\pi t/2}\Bigl(\frac f{2\pi}\Bigr)^{1/2+it}\Ga(\tfrac12+it)=\sqrt f\,\e^{i(t\log t-t+t\log(f/2\pi))}(1+O(t^{-1})).
\]
The phases cancel exactly: $t\log t-t+t\log\frac f{2\pi}+t-t\log\frac{tf}{2\pi}=0$. The product is therefore
$2\pi\GQ(f,-1)f^{-3}\cdot f\,\e^{i\pi/4}(1+O(T/K^2+T^{-1/3}))$, and $|c_Q|$ follows from Lemma~\ref{lem:Gsq}. Finally
$\int_{-H}^H\e^{-(u/H_0)^2}du=\sqrt\pi H_0(1+O(\e^{-T^{2\eps}}))$.
\end{proof}

So for even $m$ the upper bound for $|I|$ contains the term $\sqrt\pi c_QH_0$, which is of the same order as the lower bound
$|I|\gg H_0$ of \cite[(2.2)]{BSS}. No bound for the remaining terms can produce a contradiction.

\begin{remark}
Other twists do not help. If $h/k\ne1/f$, the condition $|x_0-x_*|\le K/2$ forces $1/(fk)\le|h/k-1/f|\ll K/T$, i.e.\ $k\gg T/K$.
Then $N_0\asymp h^2K^2/T\gg T$: the dual sum is longer than the original sum. Whether $\zeta_Q$ has zeros in intervals
$[T,T+T^\theta]$ with $\theta<\frac12$ for every even $m$ is open; for $m=2$ the next proposition gives bounded gaps.
\end{remark}

\begin{proposition}\label{prop:x24y2}
For $Q=x^2+4y^2$ ($m=2$),
\[
\zeta_Q(s)=2\bigl(1-2^{-s}+2^{1-2s}\bigr)\zeta(s)L(s,\chi_{-4}),
\]
and on $\sigma=\frac12$ the first factor equals $\e^{-i\theta}(2\cos\theta-2^{-1/2})$ with $\theta=t\log2$. In particular $\zeta_Q$ has zeros
$\frac12+i\gamma$ in every interval of length $2\pi/\log2$.
\end{proposition}
\begin{proof}
Let $r_2(n)=\#\{x^2+y^2=n\}=4a(n)$ with $a=1*\chi_{-4}$, so that $a(2n)=a(n)$. Since $x^2+4y^2\equiv x^2\pmod4$: for odd $n$ the
solutions of $x^2+y^2=n$ have exactly one even coordinate, so $r_Q(n)=r_2(n)/2$; for $n\equiv2\pmod4$, $r_Q(n)=0$; for $4\mid n$,
$x$ is even and $r_Q(n)=r_2(n/4)$. Hence
\[
\zeta_Q(s)=2(1-2^{-s})\zeta(s)L(s,\chi_{-4})+4^{1-s}\zeta(s)L(s,\chi_{-4}).
\]
On $s=\frac12+it$, with $z=2^{-s}=2^{-1/2}\e^{-i\theta}$, we get $1-z+2z^2=\e^{-i\theta}(2\cos\theta-2^{-1/2})$, which vanishes when
$\cos\theta=2^{-3/2}$.
\end{proof}

\begin{remark}
(i) In Theorem~\ref{thm:square} the dual sum was estimated trivially. A power saving for these sums, which have length
$\ll K^2/T$, would lower the exponent below $1/3$; we do not pursue this.
(ii) For $Q=x^2+y^2$ ($m=1$), $\zeta_Q=4\zeta(s)L(s,\chi_{-4})$, and gaps $T^{1/4+\eps}$ already follow from the zeros of $\zeta$
(Hardy--Littlewood). For $m=3,5,7,\dots$ the class number $h(-4m^2)$ exceeds one, and $\zeta_Q$ has no Euler product; e.g.\
$\Delta=36$, $Q\in\{x^2+9y^2,\ 2x^2+2xy+5y^2\}$.
\end{remark}

\section{Further remarks}\label{sec:remarks}

\begin{remark}[One-dimensional differencing; a correction to {\cite{BSS}}]\label{rem:BSS}
The argument of \cite[\S4]{BSS} differences only in one direction $e$, with shifts $\ell=me$, $1\le m\le M$. The diagonal is then
$N^2/M$ instead of $N^2/M^2$, and the bound \eqref{eq:Sl} for $\ell=(0,2m)$ reproduces \cite[(4.25)]{BSS}, giving
$S(N)\ll N^{3/4}K^{1/8}\log T$, i.e.\ the exponent $3/7$.

In \cite[(4.5)]{BSS} the difference $f_x(y+m)-f_x(y-m)$ is replaced by $2mF_x(y)$ with an error
$O(mN^{3/2}K^{3/2}T^{-2}+m^3K^{1/2}N^{-1/2})$ for each $x$. When these errors are inserted into \cite[(4.9)]{BSS}, the sum over
$x$ (of length $\asymp N^{1/2}$) has to be taken into account. The two corresponding terms in \cite[(4.9)]{BSS} should
therefore read
\[
MN^3K^{3/2}T^{-2}\quad\text{and}\quad NM^3K^{1/2}\qquad\text{instead of}\qquad MN^{5/2}K^{3/2}T^{-2}\quad\text{and}\quad M^3N^{1/2}K^{1/2}.
\]
(The B-process term $N^2/(M^{1/2}K^{1/4})$ in \cite[(4.9)]{BSS} does include the $x$-sum.) With $M=N^{1/2}K^{-1/4}$ the corrected
term $NM^3K^{1/2}=N^{5/2}K^{-1/4}$ exceeds $N^{3/2}K^{1/4}$ when $N>K^{1/2}$.

The remedy is not to expand: put $F_{x,m}(y)=(f_x(y+m)-f_x(y-m))/(2m)$ exactly. The cubic Taylor term changes $F''_{x,m}$ and
$G''_{m,n}$ only by a relative $O(m^2/N)=O(K^{-1/2})$ for $m\le M$, since each derivative costs a factor $N^{-1/2}$. Hence
\cite[(4.6), (4.16), (4.22)]{BSS} remain valid, the term $NM^3K^{1/2}$ disappears from \cite[(4.26)]{BSS}, and the exponent
$3/7$ of \cite{BSS} stands. The present paper uses no Taylor expansion; Lemma~\ref{lem:hess} treats the exact difference.
\end{remark}

\begin{remark}[On {\cite[\S4]{JS}}]\label{rem:JS}
In \cite[\S4]{JS} it is used that at every point one of $F_{xxx},F_{yyy}$ is $\asymp K^{1/2}N^{-1}$. By \eqref{eq:der},
$\partial_{e}^3\sqrt P$ vanishes where $\det(w,e)=0$ or $B(w,e)=0$. For $e=(1,0),(0,1)$ both vanish on the coordinate axes when the
coefficient of $xy$ in $P$ is $0$ (e.g.\ $P=x^2+y^2$). The repair: for fixed $w\ne0$ a direction $e$ is bad only if it lies on one
of the two lines $\R w$ and $\{e:B(w,e)=0\}$, so one of the three pairwise independent directions in $E$ is good. By compactness
there are finitely many cones covering $\R^2\setminus\{0\}$ on each of which a fixed $e\in E$ has $|\det(w,e)|\gg|w|$ and
$|B(w,e)|\gg|w|$. Changing variables by $G\in SL_2(\Z)$ with $Ge_2=e$ on each cone, the argument of \cite{JS} goes through and
gives $5/11$.
\end{remark}

\begin{remark}[The limits of the method]\label{rem:limits}
\emph{One-dimensional schemes.} If, in the scheme of \cite{BSS}, the $x$-sum is estimated by an exponent pair $(p,q)$ instead of
the second-derivative test, the bound for each $m$ becomes $(m\kappa)^{1/2+p}N^{(q-p)/2}$. Balancing
$NM^{1/2+p}K^{1/4+p/2}N^{(q-p)/2}$ against $N^2/M$ at $N=K$ gives $\sigma=\frac12E$ with $E=\frac{3/4-q/2}{3/2+p}$. The pair
$(\frac12,\frac12)$ gives $\sigma=\frac18$. The pair $(\frac27,\frac47)=BA^2B(0,1)$ would give $\sigma=\frac{13}{100}$, $\theta=\frac{37}{87}$, but
requires non-vanishing of higher derivatives of the phase $G_{m,n}$ that we have not verified. Even the conjectural pair
$(0,\frac12)$ gives only $\sigma=\frac16$. So no one-dimensional estimate of the $x$-sum beats $2/5$ in that scheme.

\emph{Why $1/3$ is out of reach for non-square $\Delta$.} $\theta=\frac13$ needs $\sigma=\frac14$, a saving $K^{1/4}$ for a
two-dimensional sum at its self-dual point. In Jutila's result for level one cusp forms \cite{J1}, and in
Theorem~\ref{thm:square}, the exponent $1/3$ comes instead from a rational critical ratio with bounded denominator, which
makes the dual sum short. For non-square $\Delta$, and always for indefinite forms, the critical ratio is a quadratic
irrational. Then $k\gg\sqrt{T/K}$ by Lagrange's theorem, and $N_0\asymp K$; taking $k$ larger only lengthens the dual sum.

\emph{Spectral reformulation.} For $Q=[a,b,c]$ positive definite and $z_Q=(-b+i\sqrt\Delta)/(2a)$ one has $Q(m,-n)=a|m+nz_Q|^2$,
hence
\[
\zeta_Q(s)=2\zeta(2s)\Bigl(\frac2{\sqrt\Delta}\Bigr)^sE(z_Q,s),\qquad E(z,s)=\frac12\sum_{(c,d)=1}\frac{(\Im z)^s}{|cz+d|^{2s}},
\]
and the zeros considered here are those of $t\mapsto E(z_Q,\frac12+it)$ at a fixed Heegner point. For indefinite $Q$, unfolding
shows that $\zQ$ is, up to $\zeta(2s)$ and elementary factors, the integral of $E(\cdot,s)$ over the closed geodesic of $Q$
(cf.\ \cite{Z}). However, the Hardy--Littlewood integral becomes a wave packet $\int h(t)E(z_Q,\frac12+it)\,dt$, i.e.\ an
incomplete Eisenstein series at $z_Q$, which is the critical Dirichlet polynomial again. The cuspidal decomposition of
$\zeta_Q$ is finite (Hecke $L$-functions of class group characters). So, for a fixed form, there is no family over which
spectral averaging could help. For $Q=x^2+4y^2$, i.e.\ $z_Q=2i$, the Hecke relation $(T_2E(\cdot,s))(i)=(2^s+2^{1-s})E(i,s)$ gives
$E(2i,s)=\frac12(2^s+2^{1-s}-1)E(i,s)$, which is equivalent to Proposition~\ref{prop:x24y2}.

\emph{Beyond Weyl.} The saving $\sigma=\frac16$ is of Weyl strength. In the $t$-aspect the Weyl barrier for $\mathrm{GL}(2)$ has been
crossed only recently, for level one forms \cite{HMSS}, and for long $\mathrm{GL}(2)$ exponential sums in \cite{HMQ}. Whether these
methods apply to the short sums arising here, with weight one dihedral coefficients, is open.
\end{remark}

\begin{remark}[Narrow ideal classes]
By Remark~\ref{rem:narrow}, for a narrow class $C^+$ with $C^+\ne C^+[(\sqrt\Delta)]$ the zeta-function $\zeta(s,C^+)=Z_Q(s)$ is not a
constant multiple of a self-dual function. Theorem~\ref{thm:B} applies to its even and odd parts separately. Whether
$\zeta(s,C^+)$ itself has infinitely many zeros on the critical line seems to be a different problem.
\end{remark}

\begin{remark}[Numerical checks]
The identities on which the transformation rests were also checked numerically: the Gauss sum formulas of
Lemma~\ref{lem:gauss} for eight definite and indefinite forms, all $k\le20$ and all $h,w$; the unified Voronoi formula
(Corollary~\ref{cor:unified}) for several forms and moduli, to between $10^{-8}$ and $10^{-13}$; the Mellin transforms of
Lemma~\ref{lem:mellin}; the stationary phase approximation of Lemma~\ref{lem:statphase} at $T=10^7$; and Lemma~\ref{lem:Gsq} for
all primitive reduced forms with $\Delta=4m^2$, $m\le30$. The scripts are available from the authors on request.
\end{remark}

\section*{Declaration of AI assistance.} The Large language model (Claude, Anthropic) has been used as assistant for writing up the manuscript, checking, and typesetting some routine details. All new ideas in the paper were conceived without the use of any AI assistant. The authors take full responsibility for the content and correctness of the paper.


\begin{thebibliography}{99}
\bibitem{BR} S.~Baier and S.~Roy, \emph{A note on the distribution of prime ideals in real quadratic fields},
arXiv:2307.03513 (2023).
\bibitem{BSS} S.~Baier, K.~Srinivas and U.~K.~Sangale, \emph{A note on the gaps between zeros of Epstein's
zeta-functions on the critical line}, Funct. Approx. Comment. Math. \textbf{57} (2017), 235--253.
\bibitem{B} B.~C.~Berndt, \emph{The number of zeros of the Dedekind zeta-function on the critical line}, J. Number
Theory \textbf{3} (1971), 1--6.
\bibitem{BKY} V.~Blomer, R.~Khan and M.~Young, \emph{Distribution of mass of holomorphic cusp forms}, Duke Math. J.
\textbf{162} (2013), 2609--2644.
\bibitem{BH} E.~Bombieri and D.~A.~Hejhal, \emph{On the distribution of zeros of linear combinations of Euler products},
Duke Math. J. \textbf{80} (1995), 821--862.
\bibitem{DH} H.~Davenport and H.~Heilbronn, \emph{On the zeros of certain Dirichlet series}, J. London Math. Soc.
\textbf{11} (1936), 181--185 and 307--312.
\bibitem{Ep} P.~Epstein, \emph{Zur Theorie allgemeiner Zetafunctionen}, Math. Ann. \textbf{56} (1903), 615--644.
\bibitem{HL} G.~H.~Hardy and J.~E.~Littlewood, \emph{The zeros of Riemann's zeta-function on the critical line},
Math. Z. \textbf{10} (1921), 283--317.
\bibitem{Iv} A.~Ivi\'c, \emph{The Riemann Zeta-Function}, John Wiley \& Sons, New York, 1985.
\bibitem{Ka} A.~A.~Karatsuba, \emph{On the distance between consecutive zeros of the Riemann zeta function that lie on the
critical line}, Trudy Mat. Inst. Steklov. \textbf{157} (1981), 49--63 (in Russian).
\bibitem{GK} S.~W.~Graham and G.~Kolesnik, \emph{Van der Corput's Method of Exponential Sums}, London Math. Soc.
Lecture Note Ser.~126, Cambridge Univ. Press, 1991.
\bibitem{HMQ} R.~Holowinsky, R.~Munshi and Z.~Qi, \emph{Beyond the Weyl barrier for $\mathrm{GL}(2)$ exponential sums},
Adv. Math. \textbf{426} (2023), Paper No.~109099.
\bibitem{HMSS} R.~Holowinsky, R.~Munshi, P.~Sharma and J.~Streipel, \emph{Sub-Weyl bound for $\mathrm{GL}(2)$ via trivial delta},
arXiv:2503.00656 (2025).
\bibitem{J1} M.~Jutila, \emph{Transformation formulae for Dirichlet polynomials}, J. Number Theory \textbf{18}
(1984), 135--156.
\bibitem{J2} M.~Jutila, \emph{Lectures on a Method in the Theory of Exponential Sums}, Tata Institute Lectures on
Math. and Phys.~80, Springer, 1987.
\bibitem{J3} M.~Jutila, \emph{Exponential sums connected with quadratic forms}, in: Number Theory (ed.\ R.~A.~Mollin),
de Gruyter, 1990, 271--286.
\bibitem{JS} M.~Jutila and K.~Srinivas, \emph{Gaps between the zeros of Epstein's zeta-functions on the critical
line}, Bull. London Math. Soc. \textbf{37} (2005), 45--53.
\bibitem{KZ} E.~M.~K\i ral and F.~Zhou, \emph{The Voronoi formula and double Dirichlet series}, Algebra Number Theory
\textbf{10} (2016), no.~10.
\bibitem{Ma} H.~Maass, \emph{\"Uber eine neue Art von nichtanalytischen automorphen Funktionen und die Bestimmung
Dirichletscher Reihen durch Funktionalgleichungen}, Math. Ann. \textbf{121} (1949), 141--183.
\bibitem{PT} H.~S.~A.~Potter and E.~C.~Titchmarsh, \emph{The zeros of Epstein's zeta-functions}, Proc. London
Math. Soc. (2) \textbf{39} (1935), 372--384.
\bibitem{Sa} A.~Sankaranarayanan, \emph{Zeros of quadratic zeta-functions on the critical line}, Acta Arith.
\textbf{69} (1995), 21--37.
\bibitem{vdD} L.~van den Dries, \emph{Tame Topology and O-minimal Structures}, London Math. Soc. Lecture Note Ser.~248,
Cambridge Univ. Press, 1998.
\bibitem{vdDM} L.~van den Dries and C.~Miller, \emph{Geometric categories and o-minimal structures}, Duke Math. J.
\textbf{84} (1996), 497--540.
\bibitem{Wat} G.~N.~Watson, \emph{A Treatise on the Theory of Bessel Functions}, 2nd ed., Cambridge Univ. Press, 1944.
\bibitem{Z} D.~Zagier, \emph{Zetafunktionen und quadratische K\"orper}, Springer, Berlin, 1981.
\end{thebibliography}
\end{document}